\documentclass[a4paper,12pt]{article}
\usepackage{amssymb}
\usepackage{latexsym}
\usepackage{amsmath}
\usepackage[toc,page]{appendix}
\usepackage{amsthm}
\usepackage{amsfonts}
\usepackage{amssymb}
\usepackage{amsmath,amscd}
\usepackage{graphicx}
\usepackage{color}

\definecolor{babyblue}{rgb}{0.54, 0.81, 0.94}
\definecolor{bittersweet}{rgb}{1.0, 0.44, 0.37}
\definecolor{brightmaroon}{rgb}{0.76, 0.13, 0.28}
\definecolor{icterine}{rgb}{0.99, 0.97, 0.37}
\definecolor{indiagreen}{rgb}{0.07, 0.53, 0.03}
\definecolor{aquamarine}{rgb}{0.5, 1.0, 0.83}
\definecolor{coolblack}{rgb}{0.0, 0.18, 0.39}
   	\definecolor{cobalt}{rgb}{0.0, 0.28, 0.67}
	\definecolor{amber}{rgb}{1.0, 0.75, 0.0}
	\definecolor{azure(colorwheel)}{rgb}{0.0, 0.5, 1.0}
	 \definecolor{shamrockgreen}{rgb}{0.0, 0.62, 0.38}
	 \definecolor{almond}{rgb}{0.94, 0.87, 0.8}
\definecolor{arylideyellow}{rgb}{0.91, 0.84, 0.42}
\definecolor{bubblegum}{rgb}{0.99, 0.76, 0.8}
\definecolor{babypink}{rgb}{0.96, 0.76, 0.76}
 	\definecolor{ballblue}{rgb}{0.13, 0.67, 0.8}

\newtheorem{Th}{Theorem}
\newtheorem{Prop}{Proposition}

\newtheorem{Lm}{Lemma}

\newtheorem{Rm}{Remark}

\newcommand{\be}{\begin{equation}}
\newcommand{\ee}{\end{equation}}
\newcommand{\bes}{\begin{equation*}}
\newcommand{\ees}{\end{equation*}}

\def\11{1\!\!1} 
\newcommand{\K}{\mathbb{H}}
\newcommand{\R}{\mathbb{R}}
\newcommand{\N}{\mathbb{N}}
\newcommand{\C}{\mathbb{C}}
\newcommand{\Z}{\mathbb{Z}}

\usepackage[colorinlistoftodos,prependcaption,textsize=tiny]{todonotes}

\usepackage[colorlinks=true,urlcolor=blue, citecolor=red,linkcolor=blue,linktocpage,pdfpagelabels, bookmarksnumbered,bookmarksopen]{hyperref}

\newcommand{\reset}{\setcounter{equation}{0}\setcounter{Th}{0}\setcounter{Prop}{0}\setcounter{Co}{0}
\setcounter{Lm}{0}\setcounter{Rm}{0}}

\def\XXint#1#2#3{{\setbox0=\hbox{$#1{#2#3}{\int}$ }
\vcenter{\hbox{$#2#3$ }}\kern-.6\wd0}}

\def\lf{\left}
\def\rg{\right}

\def\eps{\varepsilon}
\def\ds{\displaystyle}

\def\om{\omega}
\def\p{\partial}

\AtBeginDocument{
  \let\div\relax
  \DeclareMathOperator{\div}{div}
}

\begin{document}

\title{Integrability by compensation for  Dirac Equation}

\author{ Francesca Da Lio, Tristan Rivi\`ere and Jerome Wettstein\footnote{Department of Mathematics, ETH Zentrum,
CH-8093 Z\"urich, Switzerland.}}

%\date{ }
\maketitle

{\bf Abstract :} {\it We consider the Dirac Operator acting on  the  Clifford Algebra ${C\ell}_{m}$.
We show that under critical assumptions on the potential and the spinor field the equation is subject to an integrability by compensation phenomenon and has a sub-critical behaviour below some positive energy threshold (i.e. $\epsilon-$regularity theorem). This extends in 4 space dimensions as well as in 3 dimensions
a similar result obtained previously by the two first authors in 2 D in \cite{DLR1}.  }

\medskip


\medskip


\tableofcontents

\section{Introduction}
The present paper is a new contribution to the study of  linear critical systems with special structures enjoying integrability by compensation properties.

 In \cite{Riv1}, the first author
proved the sub-criticality of local {\it a-priori} critical Sch\"odinger systems in $2$ dimensions of the form
\begin{equation}
\label{zz2}
 -\Delta u= \Omega \cdot\nabla u \,~~\mbox{in ${\mathcal{D}}^{\prime}(B^2),$}
\end{equation}
where $u=(u^1,\cdots,u^{n})\in W^{1,2}(B^2,{\R}^n)$ and $\Omega\in L^2(B^2,{\R}^2\otimes \mathfrak{so}(n))$, ($\mathfrak{so}(n)$ is the Lie algebra of antisymmetric $n\times n$ matrices).
  Systems of the form  \eqref{zz2} are  related to concentration compactness and regularity results  of  Euler-Lagrange equations of conformal invariant functionals in $2$-D, such as for instance  the   {\em harmonic map equation}.\par
 Following \cite{Riv1}, in a series of works, various critical local and non local systems with antisymmetric potentials, often related to geometric variational problems, have been singled out as  enjoying  compactness properties similar to the ones of (\ref{zz2}). Successively the following systems for the corresponding critical regimes\footnote{In the function spaces which makes them critical.} and where $\Omega$ denotes an antisymmetric potential have been proven to have subcritical 
 behaviour below a threshold of energy
  \be
 \label{ex-4d}
 \Delta^2 u=\Delta(V\cdot\nabla u)+\mbox{div}(w\, \nabla u)+\Omega\cdot\nabla u
 \ee
in \cite{LR}, 
\be
\label{ex-1}
-\Delta v=\Omega\, v
\ee
in \cite{Riv2}. 

In the  nonlocal  framework, denoting
for  $\sigma\in (0,1)$ \[
(-\Delta)^{\sigma} u(x)=PV\int_{\R}\frac{u(x)-u(y)}{ |x-y|^{1+2\sigma}}\ dy\ ,
\]
similar sub-critical behaviour have been proven to hold for systems of the form respectively
\be
\label{ant1}
(-\Delta)^{1/4} v= \Omega \cdot  v \,, 
\ee
in \cite{DLR2} , as well as
\be
\label{ant2}
 (-\Delta)^{1/2} u=\Omega\cdot d_{1/2}(u)\ ,
\ee
in \cite{MS} where $d_{1/2}$ is the  the half gradient given by
\[
d_{1/2}\varphi(x,y)=\frac{\varphi(x)-\varphi(y)}{|x-y|^{1/2}}\ ,
\]
as well as
\be
\label{self}
 (-\Delta)^{1/4} u=\int_{\R} K(x,y)\ u(y) dy\ , 
 \ee
 where $ K(x,y)=-\,K^t(y,x)$ (see \cite{DLR3}).

 \medskip

In all the above examples the antisymmetry  \eqref{ant1}, \eqref{ant2} or the anti-self adjoint duality \eqref{self}  of  the potential appearing in the equation
are responsible for the regularity of the solutions or for the stability  under weak convergence as  in the original work   \cite{Riv1}. Recently in \cite{DLR1} the first and the second authors have discovered new integrability by compensation phenomena for linear systems in $2$-D where the antisymmetry is not directly involved. They are systems of the form
\begin{equation}
	\label{divergencepde}
		\operatorname{div} \left( S\, \nabla u \right) = 0 ~~~\mbox{in ~~${\mathcal{D}}^{\prime}(\R^2),$}
	\end{equation}
	where $u \in L^{2}(\C), \mathbb{R}^{n})$ and
	$S \in W^{1,2}(\C, Sym(n))$ where  $Sym(n)$ denotes the set of symmetric $n \times n$-matrices over $\mathbb{R}$ and where the crucial {\it involution} assumption is made
\be
\label{involution}
S^2=Id_n\ .
\ee 
In the case of $2$-D codomains ($n=2$) the resolution of \eqref{divergencepde} required  a different formulation of the equation   in the form
\begin{equation}\label{systconj}
	\partial_{z} \frak{f} = \Omega\cdot \overline{\frak{f}}~~\mbox{in ${\mathcal{D}}^{\prime}(\C),$}
\end{equation}
 where    $\Omega\in L^2(\C,\mathfrak{so}(2)\otimes{\C})$ is given by
 \begin{equation}\label{omegaintr}
\Omega=\left(\begin{array}{cc} 0 & \beta\\ -\beta& 0\end{array}\right),\end{equation}
for some $\beta\in L^2(\C,\C)$
  and $\frak{f} \in L^{2}(\C, \mathbb{C}^{2})$, $\bar{\frak{f}}=(\bar{\frak{f}^1},\bar{\frak{f}^2})$, and $\bar{\frak{f}^i}$ is    the complex  conjugate of $\frak{f}^i$, (see Proposition III.2 in \cite{DLR1}).
 We observe that in this context  the Lie Algebra $\mathfrak{so}(2)\otimes{\C}$ does \underbar{not} generate a \underbar{compact} Lie group. This differs completely from all the previously mentioned results above where the compactness of the underlying Lie Group, $SO(n)$, was the crucial assumption allowing the construction of suitable gauge transformations {\it \`a la Uhlenbeck}  \cite{Uh}  in order to  ``absorb'' the potential in the left-hand-side of the system.\par
The main result of \cite{DLR1} leading to the regularity of solutions  to (\ref{divergencepde}) for $n=2$ is the following theorem.
 \begin{Th}[Theorem III.7 in \cite{DLR1}]\label{regul}
 Let    $\beta\in L^2(\C,\C)$ with 
 \[
 \p_{x_1}\beta_2-\p_{x_2}\beta_1=0\ .
 \]   
 Let $\frak{f}\in L^{2}(\C,\C\times\C)$ be a solution of
 \be
 \label{main-equ-2}
 \partial_z\frak{f}= \left(
 \begin{array}{cc}
 0 &\beta\\[3mm]
 -\beta & 0
 \end{array}
 \right)\ \overline{\frak{f}}
 \ee
  Then  $\frak{f}\in L^ \frak{q}_{loc}(\C)$ for all $q<\infty$.   \hfill $\Box$
 \end{Th}
 We observe that actually the system \eqref{systconj} is critical  in the sense that if we start with a $L^{2}(\C)$ solution $\frak{f}$ then from the fact that $\partial_L\frak{f}\in L^1$ we 
get that $\frak{f}\in L^{2,\infty}(\C)$ namely we return almost to the starting point.  The   new integrability by compensation results discovered in \cite{DLR1}   are related to Wente's inequality for $2$-D Jacobians.\par

\medskip

The purpose of the present work is to extend the integrability by compensation result given by Theorem~\ref{regul} to higher dimensions. For this purpose we need to recall the fundamental notions related to Clifford Algebras:

   For every $m\ge 0$, we denote by $  {C\ell}_{m}$ the universal Clifford algebra on $\R^m$ (sometimes also denoted  ${C\ell}(0,m)$).  $  {C\ell}_{m}$ is a real associative algebra with identity containing linearly a copy of $\R^m$, such that   for any orthonormal basis $(e_1,\ldots,e_{m})$ of $\R^m$, it holds $$e_ie_j+e_je_i=-\delta_{ij},$$ for $1\le i,j\le m-1$ and the reduced  products 
$e_I=e_{i_1}\cdots e_{i_k},$ $1\le i_1\le \cdots\le e_k\le m$ and $e_{0}=1$ are a basis for $  {C\ell}_{m}$.
\footnote{If $m=0,1,2$ then $  {C\ell}_{0}\simeq\R$,  $ {C\ell}_{1}\simeq\C$ and  ${C\ell}_{2}\simeq \K$ respectively,  where \[\K:=\{\displaystyle a+b \cdot   {i} +c \cdot   {j} +d \cdot   {k},~~(a,b,c,d)\in\R^4\}, \] is the algebra of quaternions.

 ${C\ell}_{3}$ is  a real 8 dimensional space with a basis given by the following {\it paravectors}
\[
\lf\{
\begin{array}{l}
\ds e_0 \quad \mbox{ \it Scalar}\\[3mm]
\ds e_1, \, e_2, \, e_3 \quad \mbox{ \it Vectors}\\[3mm]
\ds e_1e_2, \ e_2e_3,\ e_3e_1 \quad \mbox{ \it Bivectors}\\[3mm]
\ds e_1e_2e_3 \quad \mbox{ \it Trivector}
\end{array}
\rg.
\]
  }
Any  $f\in   {C\ell}_{m}$ can be decomposed as follows:
\[
f=\sum_{I}f_{I}\,e_{I},\quad\mbox{ where }\ e_I=e_{i_0}e_{i_1}\cdots e_{i_k}\ , \ I=\{i_0,\ldots,i_k\}\quad, \quad 0\le i_1\le \cdots\le  i_k\le m \ .
\]
Let  $\sigma\colon{C\ell}_{m}\to {C\ell}_{m}$ be  the unique involutive automorphism such that $\sigma(e_i)=-e_i$ for every $i=1,\ldots,m$  and $\sigma|_{\R}=Id$. it is called the {\em principal automorphism} on ${C\ell}_{m}$ in mathematics and {\em grade involution } or {\em grade automorphism}  in physics\footnote{https://en.wikipedia.org/wiki/Paravector  }  . For $\frak{f}\in   {C\ell}_{m}$ we also denote
 \[
\hat{\frak{f}}:=\sigma(\frak{f})\ .
\]
Observe for instance  that by definition  \[
\hat{e_0}=e_0\ , \ \hat{e_i}=-e_i\  \ \widehat{e_ie_j}=\hat{e_i}\hat{e_j}=e_ie_j\ \cdots
\]  We point out that the  principal automorphism $\sigma$ is the only involution which is compatible with the Clifford Algebra structure\footnote{In the case $m=1$ then the principal automorphism coincides with the complex conjugation: $\hat{f}=\bar f$. While in the case of $m=2$ with ${C\ell}_2\simeq\K$ the automorphism $\sigma$ \underbar{does not} coincide with the other involution on $\K$ given by the conjugation operation on quaternions :
\[
\bar{1} =1\ ,\ \bar{i}=-i\ ,\ \bar{j}=-j\ \mbox{ and }\ \bar{k}=-k
\]
while
\[
\hat{ 1}=1=e_0\ ,\ \hat{ i}=-i=e_1\ ,\ \hat{ j}=-j=e_2\ ,\ \hat{ k}=k=e_1\,e_2\ 
\]}.  We will refer for instance to  \cite{garling} and \cite{hamilton} for a presentation of Clifford Algebras.

\medskip

Finally  for $\frak{f}\colon\R^{m}\to ({C\ell}_{m})^2$ we consider  the 
  Dirac operator  $\partial_L \frak{f}$   defined by  
       \begin{equation}\label{soperatorCRF}
	\partial_{L} \frak{f} = e_{0}\cdot\partial_{x_{0}} \frak{f} - e_{1} \cdot \partial_{x_{1}} \frak{f} - \ldots - e_{m-1} \cdot \partial_{x_{m-1}} \frak{f}\ ,
\end{equation}
Our main result in the present work is the following integrability by compensation theorem which  is the 3 and 4 dimensional counterpart of theorem~\ref{regul}
 \begin{Th}\label{regul-quatre}
 Let $m=3,4$.   $\beta=(\beta_0,\cdots\beta_{m-1})\in W^{1,m/2}_{loc}({\R}^m,\operatorname{span}_{\R}\{e_0,\cdots,e_{m-1}\})$ with 
 \be
 \label{curl-0}
 \forall \ i,j=1\cdots m-1\quad\quad \p_{x_i}\beta_j-\p_{x_j}\beta_i=0\ .
 \ee   
 Let $\frak{f}\in L^{m/m-1}({\R}^m,{C\ell}_{m-1}\times {C\ell}_{m-1})$ be a solution of
 \be
 \label{main-equation}
 \partial_L\frak{f}= \left(
 \begin{array}{cc}
 0 &\beta\\[3mm]
 -\beta & 0
 \end{array}
 \right)\ \hat{\frak{f}}
 \ee
  Then  $\frak{f}\in L^ {q}_{loc}(\R^m,{C\ell}_{m-1}\times {C\ell}_{m-1})$ for all $q<\infty$.   \hfill $\Box$
 \end{Th}
\begin{Rm}
Observe that the system (\ref{main-equation}) is critical in the sense that the r-h-s is a-priori only\footnote{Indeed we have $W^{1,m/2}_{loc}({\R}^m)\hookrightarrow L^{m}({\R}^m)$ } in $L^1_{loc}$ which is preventing a direct use of Calderon Zygmund theory. Any direct attempt to bootstrap is blocked by the fact that  $\partial_L^{-1}L^1_{loc}\hookrightarrow L^{m/m-1,\infty}_{loc}({\R}^m)$. Which means that a-priori integrability information
on $\frak{f}$ is lost from the first iteration on. It is only because of its very peculiar structure that, thanks to some ``hidden'' compensation, a gain of integrability and local compactness holds. In fact a quantitative
version of the theorem~\ref{regul-quatre} can bee formulated in the form of an $\epsilon$-regularity. \hfill$\Box$
\end{Rm}
\begin{Rm}
Some gain of integrability still holds when instead of assuming (\ref{curl-0}) one assumes that $\p_{x_i}\beta_j-\p_{x_j}\beta_i\in L^p_{loc}({\R}^m)$ ($m=3,4$) for some $p>2$.\hfill $\Box$
\end{Rm}
\begin{Rm}
It would be interesting to study the possibility whether or not theorem~\ref{regul-quatre} continues to hold if instead of assuming $\beta$ to be in $W^{1,m/2}_{loc}$ one would make the 
milder hypothesis $\beta\in L^m_{loc}$. In fact, we are proving theorem~\ref{regul-quatre} under the assumption that $\beta$ belongs to the Lorentz space $L^{(m,2)}_{loc}$ in which $W^{1,m/2}_{loc}$ embeds in $m$ dimensions for $m=3,4$(see \cite{Riv-notes}). \hfill$\Box$\end{Rm}
\begin{Rm}
The investigations made by the authors is leading them to the conclusion that the theorem does not generalise to arbitrary $m$ in a straightforward way and the proofs given below for the cases $m=3,4$ is very much
``dimension depending''. Some results have been obtained by the third author in dimensions $m\le 8$ in \cite{Wett}. 
\end{Rm}
Similarly to the 2-dimensional case the resolution of Theorem~\ref{regul-quatre} for $m=4$ for instance goes through the canonical inclusion of ${C\ell}_{3}$ into ${C\ell}_{4}$ (i.e. ${C\ell}_{4}\simeq{C\ell}_{3}\oplus {C\ell}_{3}\, e_4$) and the introduction of
the new variable $\frak{g}=\frak{f}^1+\frak{f}^2 e_4\in {C\ell}_{4}$. The equation satisfied by $\frak{g}$ is then\footnote{In this form the equation identifies to the {\bf covariant Dirac} equation commonly written as follows
\be
\label{cov-dirac}
\sum_{\mu=0}^3\gamma_\mu (\p_{x_\mu}-\mathbf{i} A_\mu)\psi=0
\ee
where $\gamma_0$ is the $2\times 2$ identity matrix, $\gamma_\mu=-e_\mu$ , $\mathbf{i}=e_4$, the connection components are given by $A_\mu=\beta_\mu$ and the group representing
on the spinor space ${C\ell}_4\simeq{C\ell}_{3}\oplus {C\ell}_{3}\, e_4$ is the abelian group $\exp(e_4{\R})$. With this identification at hand one could then imagine that, for instance assuming $\beta_0=0$ and $\p_{x_0}\beta=0$ the  flatness of the connection $A$ implied by (\ref{curl-0}) would make the absorption of the r-h-s of (\ref{syst-g}) trivial by multiplying on the left (\ref{cov-dirac}) by $\exp(e_4\varphi)$ where $d\varphi=A$. However we have
\[
\forall \ l=1,2,3\quad\quad\exp(e_4\,\varphi) e_l=e_l\, \exp(-e_4\,\varphi)
\]
and this multiplication would then give
\[
\p_{x_0}(\exp(e_4\,\varphi) \frak{g})-\sum_{i=1}^3e_l\ \p_{x_i}\lf(\exp(-e_4\,\varphi) \frak{g})\rg)=0\ ,
\]
which is not easily invertible neither unless in the very particular case where $\frak{g}$ is known to be independent of $x_0$ (that we are not assuming a-priori).} 
\be
\label{syst-g}
\partial_L\frak{g}=-(\beta e_4)\cdot\frak{g}~~\mbox{in ${\mathcal{D}}^{\prime}(\R^4).$} 
\ee
The ``absorption'' of the right hand side of this equation by the left-hand-side will be achieved through the construction {\it \`a la Uhlenbeck} of a  {\it Coulomb type Gauge} in the Lie group whose Lie algebra is given by
\[
{\mathcal{E}}_4=\{e_4,e_1e_4,e_2e_4,e_3e_4,e_1,e_2,e_3,e_1e_2,e_1e_3,e_2e_3\}\ .
\]
This Lie group happens to be isomorphic to $Spin(5)$ (see Appendix for a presentation of $Spin(m)$) and is hence \underbar{compact} which is crucial for the gain of integrability similarly to the seminal work \cite{Riv1}. 

We also would like to stress that the linearized natural {\it Coulomb type condition}  in the present framework  is given by the {\bf Lorenz  gauge equation} for an electric potential 
$\varphi$  and   magnetic vector potential $A$ (see (\ref{bt8}) and (\ref{bt8-aa})) :
\be
\label{maxwel}
\lf\{
\begin{array}{rcl}
{\mathbf{curl}}(A)&=&{\mathbf{B}}\\
-\partial_t A-\nabla \varphi&=&{\mathbf{E}}\\
\partial_t \varphi+\mathbf{div}(A)&=&\beta_0
\end{array}
\rg.
\ee
where ${\mathbf{B}},$ and ${\mathbf{E}}$ represent respectively the  {\bf electric} and the {\bf magnetic fields} and are taken in our case to be ${\mathbf B}=0$ and ${\mathbf E}=(\beta_1,\beta_2,\beta_3)$ while assuming (\ref{curl-0}) and $x_0$ is the time variable.

\medskip

Finally it is legitimate to ask if the resolution of theorem~\ref{regul-quatre} leads to any new result regarding solutions in $L^{m/m-1}$ of real Elliptic Systems in Divergence form (\ref{divergencepde}) involving {\it critical chirality operator} $S\in W^{2,m/2}_{loc}({\R}^m,{Sym}(m))$ with $S^2=Id_m$. We have not been able to establish this connection so far\footnote{Actually in dimension $m=3$  $L^{3/2}$ solutions of
\begin{equation}\label{maineqbis-1}
\div(S\,\nabla u)=0~~\mbox{in ${\mathcal{D}}^{\prime}(\R^3)$}
\end{equation}
happen to be rather related to a solution of a system of the type:
\begin{equation}\label{systf-1}
\partial_L \frak{f} = \left(\begin{array}{cc} 0 & \beta\\ -\beta& 0\end{array}\right)\  \overline{{\frak{f}}}
\end{equation}
where  $\frak{f}=(\frak{f}^1,\frak{f}^2)\in \K\times \K$ and $\overline{{\frak{f}^{i}}}$ denotes the conjugate of $\frak{f}$ in ${\K}$  which does not coincide with  the involution $\sigma$ in ${C\ell}_2\simeq\K$ given by the principal automorphism  we are considering in (\ref{main-equation}).
Up to our knowledge,  the question of  a possible higher integrability for $L^{3/2}$ solutions of systems of the form \eqref{systf-1} for $\beta\in L^3$ or even in $W^{1,3/2}$ as well as the question of the possible higher integrability for $L^{3/2}$ solutions of \eqref{maineqbis-1} 
with $W^{2,3/2}_{loc}$ Involution operator $S$ that would naturally extend to 3D the theorems in 2 dimension obtained in \cite{DLR1} are still open.}.
%\newpage

\section{Bootstraptest for \eqref{syst-g} in $4$-D}\label{4Dim}
\reset
  
 We first start with the dimension $m=4$ that looks more natural to us. We consider  systems of the form
\begin{eqnarray}\label{maineq4d}
\partial_L\frak{f}&=&\beta e_4 \cdot  \frak{f}\,.
\end{eqnarray}
where  $\frak{f}=\frak{f}^1+\frak{f}^2 \, { {e_4}}  $, 
and $\frak{f}^1,\frak{f}^2\colon \R^4 \to   {C\ell}_{3}$. The function $\frak{f}$ assumes values in the Clifford algebra ${C\ell}_4.$\par
\par
In this context $\partial_L $ and $\partial_R$ denote respectively  the 
the left and right Dirac operator in $\R^4$  defined by 
\begin{eqnarray}
\partial_L \frak{f}&:=& \partial_{x_0}\frak{f}-\sum_{i=1}^3e_i \partial_{x_i}\frak{f}\label{deltaL4}\\
\partial_R \frak{f}&:=& \partial_{x_0}\frak{f}-\sum_{i=1}^3\ \partial_{x_i}\frak{f}e_i \label{deltaR4}.
\end{eqnarray}
\par
 
The first main result is the following Theorem which corresponds to a bootstrap test for the equation \eqref{maineq4d}:
\begin{Th}\label{mainth4D} There exists $\varepsilon_0>0$ such that for every $\beta\in L^{(4,2)}(\R^4, V_3 )$ satisfying $\| \beta\|_{L^{(4,2)}(\R^4)}\le \varepsilon_0$ as well as 
$$\forall i,j  \quad \partial_{x_{i}} \beta^{j} - \partial_{x_{j}} \beta^{i} = 0\ ,$$
and where $V_3={\operatorname{span}}_{\R}\{e_0,e_1,e_2,e_3\}$ and every  $\frak{f} \in L^{4/3}(\R^4,{C\ell}_4)$ solving
 \begin{equation}\label{maineq6}
\partial_L \frak{f}=\beta  e_4 \cdot \frak{f}\ ,
\end{equation}
we have ${\frak{f}}\equiv 0.$ \hfill $\Box$
\end{Th}
In order to prove Theorem \ref{mainth4D}, we first perform the construction of a suitable gauge.
This relies  on the use of certain projections to render the emerging gauge equations elliptic and therefore enabling direct existence and regularity arguments. The arguments are given in the following subsections and we will make use of a new compensation phenomenon linked to the appearance of Maxwell-type equations for our changes of gauge.    
\subsection{Construction of a Gauge}
  In order to employ an absorption argument by a change of gauge, we consider the compact Lie algebra generated  by  $\{e_4,e_1e_4,e_2e_4,e_3e_4\}$.  Such an algebra is isomorphic to $\frak{spin}(5)$  and it is given by:
$$
 {{E}}=\operatorname{span}_{\R}\{e_4,e_1e_4,e_2e_4,e_3e_4,e_1,e_2,e_3,e_1e_2,e_1e_3,e_2e_3\}.$$
We may split $${{E}}=E_4\oplus E_6,$$ where
$$E_4:=\operatorname{span}_{\R}\{e_4,e_1e_4,e_2e_4,e_3e_4\},$$
and
$${E}_6:=\operatorname{span}_{\R}\{e_1,e_2,e_3,e_1e_2,e_1e_3,e_2e_3\}.$$
Note that ${E}_6$ is also a compact Lie algebra of dimension $6$ isomorphic to $\frak{spin}(4).$
Let us introduce the following projections:
\begin{eqnarray}
&&\Pi_{6}\colon {C\ell}_4\to {E_6}\label{pi6}\\
 &&\Pi_{4}\colon {C\ell}_4\to {E}_4\label{pi4}\\
 &&\Pi_{3}\colon {C\ell}_4\to E_3:=\operatorname{span}_{\R}\{e_2e_{3}e_4, e_3e_1e_4, e_1e_2e_4\}\label{pi4}\\
 &&{\mathcal{P}}\colon E_3\to E_4,~~e_{k+1}e_{k-1}e_4\mapsto e_ke_4, ~~k=1,2,3.\label{p3}
  \end{eqnarray}
  In the projection ${\mathcal{P}}$, we use the indexing in $\Z/3\Z$. This means for example that we identify $4$ with $1$ in \eqref{p3}.
\par

 We multiply both sides of the equation \eqref{maineq6} from the left by $\frak{q}$   belonging to the compact Lie group    corresponding to the Lie algebra $E$. Such a Lie group is isomorphic to $Spin(5)\simeq Sp(2)$\footnote{see Thm. 9.11.(iii.) in \cite{garling}. The symplectic group, $Sp(n)$ is
the subgroup of $Gl(n,\K)$, the invertible $n\times n$ quaternionic matrices  $A$ satisfying $\bar{A}^t A=A \bar{A}^t=\11$.}.
 We obtain 
 \begin{equation}\label{qf}
 \frak{q}\partial_L \frak{f}=\partial_{x_0}(\frak{q} f)-(\partial_{x_0}\frak{q})\frak{f}-\sum_{i=1}^{3}\partial_{x_i}(\frak{q} e_i\frak{f})+\sum_{i=1}^{3}\partial_{x_i}\frak{q} e_i \frak{f}
 \end{equation}
 We define the following expression:
 $$
 {\mathcal{D}}(\frak{q}):=\frak{q}^{-1}\partial_{x_0}\frak{q} - \sum_{i=1}^{3}\frak{q}^{-1}\partial_{x_i}\frak{q} e_i = \frak{q}^{-1} \partial_R \frak{q}.$$
Let us observe that
 \begin{eqnarray}
\beta e_4=\beta^0\ e_4-\sum_{i=1}^3 \beta^i\  e_i e_4~~&\in& E_4\label{decalpha}\\
 {\mathcal{D}}(\frak{q})&\in& E_3\oplus E_4\oplus E_6 \oplus \R\cdot e_1 e_2 e_3.\label{D}
 \end{eqnarray}
By combining \eqref{maineq6}, \eqref{qf} and $\eqref{D}$ we get
 \begin{eqnarray}\label{neweq}
 \partial_{x_0}(\frak{q} \frak{f})-\sum_{i=1}^{3}\partial_{x_i}(\frak{q} e_i\frak{f})&=&\frak{q}(\beta\,e_4+{\mathcal{D}}(\frak{q})) \frak{f}.
 \end{eqnarray}
 We notice that \eqref{neweq} is a system of $15$ equations in $10$ unknowns, if we split the PDE according to the basis in ${C\ell}_{4}$. If we  try to directly solve:
 $$\beta\,e_4+{\mathcal{D}}(\frak{q}) = 0,$$
 this will have  therefore little to no chance of success. Instead, let us try and approximately solve this equation.\par
  More precisely, our main goal is to find $\frak{q}\in \dot W^{1,4}(\R^4, Spin(5))$\footnote{$\dot W^{1,4}(\R^4, Spin(5))$ denotes the space of functions $u\in L^1_{loc}(\R^4, Spin(5))$ such that $\nabla u\in L^4(\R^4)$} such that
   ${\mathcal{D}}(\frak{q})=-\beta e_4+{\mathcal{V}}(x)$, where $\mathcal{V}$ is a more regular potential than $\beta$, namely  ${\mathcal{V}}\in L^{(4,1)}(\R^4).$\\
   
      To achieve this, we  introduce the following operator:
   \begin{eqnarray}\label{estaux}
&& {\mathbf{N}}\colon \dot W^{1,4}(\R^4 ,Spin(5))\to  W^{-1,4}(\R^4 ,E_6)\times  L^4(\R^4,E_4)\times L^4(\R^4,E_3)\\[5mm]
  &&  ~~~\frak{q}\mapsto\left(\Pi_6\left(\sum_{i=0}^3(\partial_{x_i}((\frak{q} )^{-1}\partial_{x_i} \frak{q} ))\right),\Pi_{4}({\mathcal{D}}(\frak{q} )),\Pi_{3}({\mathcal{D}}(\frak{q} ))\right)\nonumber
  \end{eqnarray}
  We observe that ${\mathbf{N}} $ is an operator from $Spin(5)$-valued maps, i.e. functions whose range has dimension $10$, to $E_6 \oplus E_4 \oplus E_3$-valued functions, namely functions taking values in a space of dimension $13$. Therefore, there is no hope of proving that it is an isomorphism. This suggests that we have to further reduce dimensionality to arrive at an operator which takes all values sufficiently close to $0$.\\
  
   Indeed, we would like to prove the following result (which, as we stated before, is a-priori impossible in the generality presented):
   \begin{Lm}\label{inj}
   There exists $\varepsilon_0>0$ such that for every $\beta\in L^4(\R^4,V_3)$     and $\|\beta\|_{L^4(\R^4)}\le \varepsilon_0,$ 
    there exists $\frak{q}\in \dot W^{1,4}(\R^4 ,Spin(5))$ such that 
   \begin{eqnarray}
   && {\mathbf{N}}(\frak{q})=(0,\beta e_4,0).\label{Nq}
   \end{eqnarray}
   together with an estimate:
   \begin{equation}\label{Nq2}
   \|\nabla\frak{q}\|_{L^4}\lesssim \|\beta\|_{L^4(\R^4)}.
   \end{equation}
 \end{Lm}
 
   Unfortunately, this strong form of a gauge is not possible. However, if we allow for an error term of slightly better integrability, which will suffice for the regularity result we are trying to establish, we can actually achieve a suitable change of gauge by using a slightly weaker gauge operator.
% {\color{red} Up to an additional potential term with improved integrability}
 \par
 \medskip
 In order to prove a weaker analogue of Lemma \ref{inj}, we first consider a different nonlinear operator:
\begin{eqnarray}\label{gaugeop}
&& {\mathcal{N}} \colon \dot W^{1,4}(\R^4 ,Spin(5))\to  W^{-1,4}(\R^4 ,E_6)\times   L^4(\R^4,E_4) \\[5mm]
&&~~~~~ \frak{q}\mapsto\left(\Pi_6\left(\sum_{i=0}^3(\partial_{x_i}((\frak{q})^{-1}\partial_{x_i} (\frak{q}))\right), (\Pi_{4}+\mathcal{P})({\mathcal{D}}(\frak{q}))\right) \nonumber
\end{eqnarray}
   Observe that
   \begin{eqnarray}
    (\Pi_{4}+\mathcal{P})({\mathcal{D}}(\frak{q}))&=& \Pi_{e_4}({\mathcal{D}}(\frak{q}))+\sum_{i=1}^3
( \Pi_{e_ie_4}+\Pi_{\mathcal{P}({e_{i+1},e_{i-1}e_4})})({\mathcal{D}}(\frak{q}))
 \end{eqnarray}
 We shall construct a gauge for the nonlinear operator in \eqref{gaugeop}:
\begin{Lm}\label{lm-invertN-a}
There are constants  $\varepsilon_0>0$ and $C>0$ such that for any choice of $\om\in W^{-1,4}(\R^4 ,E_6)$ and $\frak{g}\in   L^{4}(\R^4,E_4)$ satisfying
\begin{equation}\label{smallass}\|\omega\|_{W^{-1,4}}\le \varepsilon_0,~~\|\frak{g}\|_{L^{4}}\le \varepsilon_0
\end{equation}
there exists $\frak{q}\in \dot  W^{1,4}(\R^4,Spin(5))$ such that  
\begin{equation}\label{4D:decomp}{\mathcal N}({\frak{q}})=(\omega,\frak{g})\end{equation}
and
\begin{equation}
\|\nabla \frak{q}\|_{L^4}\le C(\|\omega\|_{W^{-1,4}}+\|\frak{g}\|_{L^{4}})\,.
\end{equation}
\end{Lm}
 
 In order to prove Lemma~\ref{lm-invertN-a}, we shall need to introduce some notation and establish a few intermediate results: As in \cite{DLR2,DLS2}, by an approximation argument, it suffices to prove Lemma~\ref{lm-invertN-a} under the assumption that $\om$ and $\frak{g}$ are slightly more regular. More precisely, we first prove Lemma \ref{lm-invertN-a} in the case $\om\in ( W^{-1,p}\cap  W^{-1,p'})({\R^4},E_6)$ and $\frak{g}\in   (L^p\cap L^{p^{\prime}})({\R^4},E_4)$ for some fixed $4<p$ and its H\"older-dual $p'=\frac{p}{p-1}$: \\

  For $\varepsilon > 0$, let us now introduce:
\be
\label{Uepsilon}
     {\mathcal{U}}_\eps  := \left \{
     \begin{array}{c}
     (\omega,\frak{g}) \in  (W^{-1,p}\cap   W^{-1,p'})({\R^4},E_6)\times (L^p\cap L^{p^{\prime}})({\R^4},E_4)\\[5mm]
   \|\omega\|_{W^{-1,4}} +\|\frak{g}\|_{L^4}\le \varepsilon
  \end{array}
  \right \}
  \ee

\noindent
For constants $\varepsilon, \Theta > 0$, let $\mathcal{V}_{\varepsilon,\Theta} \subseteq {\mathcal{U}}_\eps$ denote the set of pairs $(\omega, \frak{g})$ for which we have a decomposition as in \eqref{4D:decomp} and which are satisfying the following estimates:
 \be
\label{4D:eq:gauge:2est}
  \|\nabla  \frak{q}\|_{L^4} \leq \Theta  ( \|\omega\|_{W^{-1,4}}+\|\frak{g}\|_{L^4})\, 
 \ee 
 \be
  \label{4D:eq:gauge:pest}
 \|\nabla \frak{q}\|_{p} \leq \Theta  ( \|\omega\|_{W^{-1,p}}+\| \frak{g}\|_{L^p})\, ,
 \ee
 \be
 \label{4D:eq:gauge:ppest} {\mathcal{U}}_\eps
     \|\nabla \frak{q}\|_{p^{\prime}} \leq \Theta  ( \|\omega\|_{W^{-1,p^\prime}}+\|\frak{g}\|_{L^{p^\prime}})\,.
\ee
That is:
\[
 \mathcal{V}_{\varepsilon,\Theta} := \left \{(\omega,\frak{g}) \in {\mathcal{U}}_\eps:\ \begin{array}{c}
                                             \mbox{there exists $\frak{q} \in (\dot{W}^{1,p} \cap \dot W^{1,p^{\prime}})(\R^4,Spin(5))$, so that  }\\[3mm]
                                           \frak{q}-\frak{I}\in  L^{4p/3p-4} (\R^4,Spin(5))\\[3mm]~~\mbox{and}~~ \eqref{4D:decomp}, \eqref{4D:eq:gauge:2est}, \eqref{4D:eq:gauge:pest},  \eqref{4D:eq:gauge:ppest}~~ \mbox{hold.}
                                              \end{array}\right \}
\]
The strategy   to prove Lemma~\ref{lm-invertN-a} follows the one introduced by K. Uhlenbeck in \cite{Uh} in order to construct Coulomb gauges in critical dimensions. 
In fact, Lemma~\ref{lm-invertN-a}  is going to be a consequence of the following proposition, which will establish our Lemma by using a standard connectedness argument on suitable spaces.
\begin{Prop}\label{4D:pr:VepseqUeps}
There exist $\Theta > 0$ and $\eps > 0$, such that $\mathcal{V}_{\eps,\Theta} = {\mathcal{U}}_\eps$. \hfill $\Box$
\end{Prop}
\noindent{\bf Proof of Proposition \ref{4D:pr:VepseqUeps}.}
Proposition~\ref{4D:pr:VepseqUeps} will follow, once we have shown the following four properties:
\begin{itemize}
 \item[(i.)] $ {\mathcal{U}}_\eps$ is  connected.
 \item[(ii.)] $\mathcal{V}_{\eps,\Theta}$ is nonempty. 
 \item[(iii.)] For any $\eps, \Theta > 0$, $\mathcal{V}_{\eps,\Theta}$ is a relatively closed subset of $ {\mathcal{U}}_\eps$.
 \item[(iv.)] There exist $\Theta > 0$ and $\eps > 0$, such that $\mathcal{V}_{\eps,\Theta}$ is a relatively open subset of $ {\mathcal{U}}_\eps$.
\end{itemize}
\par
Property (i.) is obvious, since $ {\mathcal{U}}_\eps$ is clearly starshaped with center $0$ and hence path-connected. Property (ii.) is also evident, since $(0,0)\in \mathcal{V}_{\eps,\Theta}$ follows by choosing the constant map $\frak{q} = \frak{I}$. Consequently, it remains to verify the latter two:\\

The closedness property (iii.) follows almost verbatim from those in \cite{DLR2} and \cite{DLS2}: Assume that $(\omega_n, \frak{g}_n), (\omega, \frak{g}) \in \mathcal{V}_{\varepsilon,\Theta}$ and moreover, $(\omega_n, \frak{g}_n) \to (\omega, \frak{g})$ and let $\frak{q}_n$ be as in the definition of $\mathcal{V}_{\varepsilon,\Theta}$, i.e. $\mathcal{N}(\frak{q}_{n}) = (\omega_n, \frak{g}_n)$ and satisfying \eqref{4D:eq:gauge:2est}, \eqref{4D:eq:gauge:pest},  \eqref{4D:eq:gauge:ppest}. Observe that $\nabla \frak{q}_n$ is bounded in $L^{p}$ and $L^{p^\prime}$. Therefore, we can extract weakly converging subsequences with limit $\frak{p}$. Furthermore, we may extract another subsequence of $\frak{q}_n - \frak{I}$ converging locally in $L^{q}$ for some $q < \frac{4p}{3p-4}$ we may choose, due to the $\dot{W}^{1,p^{\prime}}$-boundedness of $\frak{q}_n$. The limit $\frak{q} - \frak{I}$ satisfies $\nabla \frak{q} = \frak{p}$ and $\frak{q}$ assumes values in $Spin(5)$ a.e.. This can be seen by extracting another subsequence of $\frak{q}_n - \frak{I}$ converging a.e. pointwise and using the closedness of $Spin(5)$. Due to the weak lower semi-continuity of the norms, we immediately obtain that $\eqref{4D:eq:gauge:2est}, \eqref{4D:eq:gauge:pest}$ and $\eqref{4D:eq:gauge:ppest}$ hold. Finally, observe that, in the distributional sense, we have the convergence:
$$\Pi_6\left(\sum_{i=0}^3(\partial_{x_i}((\frak{q}_n)^{-1}\partial_{x_i} (\frak{q}_n))\right) \to \Pi_6\left(\sum_{i=0}^3(\partial_{x_i}((\frak{q})^{-1}\partial_{x_i} (\frak{q}))\right),$$
as well as
$$(\Pi_{4}+\mathcal{P})({\mathcal{D}}(\frak{q}_n)) \to (\Pi_{4}+\mathcal{P})({\mathcal{D}}(\frak{q})).$$
This shows $\mathcal{N}(\frak{q}) = (\omega, \frak{g})$ and therefore relative closedness is established. This takes care of (iii.).
\medskip

Lastly, we verify the {openness property} (iv.). For this let $\omega_0,\frak{g}_0 $ be arbitrary in $\mathcal{V}_{\eps,\Theta}$, for some $\eps, \Theta > 0$ determined later on in an appropriate manner: Let $\frak{q}_0 \in \dot{W}^{1,p} \cap \dot{W}^{1,p^{\prime}}(\R^4,\frak{spin}(5))$, $ \frak{q}_0-\frak{I} \in L^{4p/3p-4}(\R^4 )$, such that the decomposition      \eqref{4D:decomp} as well as the estimates  \eqref{4D:eq:gauge:2est}, \eqref{4D:eq:gauge:pest} and \eqref{4D:eq:gauge:ppest}  are satisfied for $\omega_0$ and $\frak{g}_0 $.

\medskip

Let us consider perturbations of $\frak{q}_0$ of the form $\frak{q}=\frak{q}_{0}e^{\frak{u}}$  
where $\frak{u}\in \dot W^{1,p} \cap \dot W^{1,p^{\prime}} \cap L^{4p/3p-4}(\R^4,\frak{spin}(5))$. Observe that the exponent $p>4$ has been chosen in particular to ensure $\frak{u}\in C^0\cap L^\infty({\R^4})$ and $\frak{q}_{0}e^{\frak{u}}-\frak{I}\in L^{\frac{4p}{3p-4}}.$ Indeed for a Schwartz function $\frak{u}$, one has
\begin{equation}
\label{contembedlppprime}
\frak{u}(x)=C\int_{\R^4}\nabla_x |x-y|^{-2}\cdot\nabla\frak{u}(y)\ dy\quad\Rightarrow\quad\|\frak{u}\|_\infty\lesssim\|\nabla_x |x-y|^{-2}\|_{L^{(4/3, \infty)}}\, \|\nabla\frak{u}\|_{L^{(4,1)}}
\end{equation}
The generalized H\"older inequality (see \cite{Gra}) yields the required estimate of the Lorentz norm:
\[
 \|\nabla\frak{u}\|_{L^{(4,1)}}\le C\,  \|\nabla\frak{u}\|^{\alpha}_{L^p}\  \|\nabla\frak{u}\|^{1-\alpha}_{L^{p^\prime}}\,.
\] 
where $4^{-1}=\alpha p^{-1}+(1-\alpha){p'}^{-1}$. The statement $\frak{u}\in L^\infty$, and thus continuity by approximation, follows due to the density of Schwartz functions in $\dot W^{1,p} \cap \dot W^{1,p^{\prime}} \cap L^{\frac{4p}{3p-4}}$. It can be easily seen, that the argument carries over to domains of arbitrary dimension $m$, if $m < p$, as the density result and the interpolation identity do not critically depend on $m = 4$ in any significant way. This observation ensures that the argument presented could be generalised up to this point to higher-dimensional domains without issues.

 \paragraph{Study of the linearized operator}
 The key idea is that we can deduce general global properties of the gauge operator by considering its differential at the element $0$. This is in line with the usual local inversion theorem, which again reduces local invertibility to a question about invertibility of the differential which is its local linearisation.\\
 
  We   compute $D {\mathcal{N}}  (\frak{q}_0)$ as follows:
\[
D {\mathcal{N}}  (\frak{q}_0)= \frac{d}{dt} {\mathcal{N}} (\frak{q}_0e^{t\frak{u}}) \Big|_{t=0}=: {\mathcal{L}}_{\frak{q}_{0}}({\frak{u}}),
\]
where $\frak{u}\in (\dot W^{1,p} \cap \dot W^{1,p^{\prime}} \cap L^{4p/3p-4})(\R^4,\frak{spin}(5))$. Let us write this in components:
%$${\mathcal{L}}_{\frak{q}_{0}}({\frak{u}})=({\mathcal{L}}^6_{\frak{q}_{0}}({\frak{u}}),{\mathcal{L}}^{4}_{\frak{q}_{0}}({\frak{u}}),{\mathcal{L}}^1_{\frak{q}_{0}}({\frak{u}}),%{\mathcal{L}}^2_{\frak{q}_{0}}({\frak{u}}),{\mathcal{L}}^3_{\frak{q}_{0}}({\frak{u}})).$$

$${\mathcal{L}}_{\frak{q}_{0}}({\frak{u}})=({\mathcal{L}}^6_{\frak{q}_{0}}({\frak{u}}),{\mathcal{L}}^{4}_{\frak{q}_{0}}({\frak{u}}))$$
where we have:
\begin{eqnarray*}
{\mathcal{L}}^6_{\frak{q}_{0}}({\frak{u}})&:=&\Pi_{6}\left[\Delta\frak{u}+
\sum_{i=0}^3\partial_{x_i}\left(\frak{q}_0^{-1}(\partial_{x_i}\frak{q}_0)\frak u-\frak{u}\frak{q}^{-1}_0\partial_{x_i}\frak{q}_0\right)\right]\\
{\mathcal{L}}^{4}_{\frak{q}_{0}}({\frak{u}})&:=&(\Pi_{4}+{\mathcal{P}})[\partial_{x_0}\frak{u}-\sum_{i=1}^{3}\partial_{x_i}\frak{u}e_i]\nonumber\\[3mm]&&+(\Pi_{4}+{\mathcal{P}})[\frak{q}_0^{-1}(\partial_{x_0}\frak{q}_0)\frak u-\frak{u}\frak{q}^{-1}_0\partial_{x_0}\frak{q}_0-\sum_{j=1}^{3}(\frak{q}_0^{-1}(\partial_{x_j}\frak{q}_0)\frak u-\frak{u}\frak{q}^{-1}_0\partial_{x_j}\frak{q}_0)e_j]
\end{eqnarray*}
and for $i=1,2,3$. First, we investigate the invertibility of ${\mathcal{L}}_{\frak{q}_0}({\frak{u}})$ in the special case $\frak{q}_0=\frak{I}.$

%\begin{eqnarray*}
%{\mathcal{L}}^{i}_{\frak{q}_{0}}({\frak{u}})&=&
%(\Pi_{e_ie_4}+\Pi_{e_{i+1}e_{i-1}e_4})[ \partial_{x_0}\frak{u}-\sum_{j=1}^{3}\partial_{x_j}\frak{u}e_j]\nonumber\\&&+
 %(\Pi_{e_ie_4}+\Pi_{e_{i+1}e_{i-1}e_4})\left[\frak{q}_0^{-1}(\partial_{x_0}\frak{q}_0)\frak u-\frak{u}\frak{q}^{-1}_0\partial_{x_0}\frak{q}_0-\sum_{j=1}^{3}(\frak{q}_0^{-1}%(\partial_{x_j}\frak{q}_0)\frak u-\frak{u}\frak{q}^{-1}_0\partial_{x_j}\frak{q}_0e_j)\right]
 %\end{eqnarray*}
 
\paragraph{Invertibility of ${\mathcal{L}}_{\frak{I}}({\frak{u}})$}
If $\frak{q}_0=\frak{I}$, we obviously have $d\frak{q} = 0$ and as a result, the operator $\mathcal{L}_{\frak{I}}(\frak{u}) = (\mathcal{L}^{6}_{\frak{I}}(\frak{u}), \mathcal{L}^{4}_{\frak{I}}(\frak{u}))$ has the following simpler form:
\begin{eqnarray}
{\mathcal{L}}^6_{\frak{I}}({\frak{u}})&=&\Pi_{6}\left[\Delta\frak{u}\right]\nonumber\\
{\mathcal{L}}^{4}_{\frak{I}}({\frak{u}})&=&(\Pi_{4}+{\mathcal{P}}) [\partial_{x_0}\frak{u}-\sum_{i=1}^{3}\partial_{x_i}\frak{u}e_i]
 \end{eqnarray}
\noindent
 The following will suffice to prove existence of solutions and regularity:

%\begin{eqnarray*}
%{\mathcal{L}}^{i}_{\frak{I}}({\frak{u}})&=&(\Pi_{e_ie_4}+\Pi_{e_{i+1}e_{i-1}e_4})[ \partial_{x_0}\frak{u}-\sum_{j=1}^{3}\partial_{x_j}\frak{u}e_j]
%\end{eqnarray*}
 
\begin{Prop}\label{invert}
 The operator ${\mathcal{L}}_{\frak{I}}({\frak{u}})$ is elliptic. \end{Prop}

We mention at this point that this will be the only point where we crucially use the dimension of the domain, as we shall observe the Riemann-Fueter operator on $\R^4$ emerging in our computations.\\

\noindent
{\bf Proof of Proposition \ref{invert}.}
We write $\frak{u}=w+v$ where $w\in  E_6$ and $v=v^0e_4+v^1e_1e_4+v^2e_2e_4+v^3e_3e_4\in E_4$.
We observe that 
\begin{eqnarray}
{\mathcal{L}}^6_{\frak{I}}({\frak{u}})&=&\Pi_{6}\left[\Delta w + \Delta v\right] = \Delta w \nonumber\\
{\mathcal{L}}^{4}_{\frak{I}}({\frak{u}})&=&(\Pi_{4}+{\mathcal{P}})\left[\partial_{x_0}v-\sum_{i=1}^{3}\partial_{x_i}ve_i \right]\nonumber
\end{eqnarray}
 We explicitly compute ${\mathcal{L}}^{4}_{\frak{I}}({\frak{u}})$:
 
 \begin{eqnarray}
 (\Pi_{4}+{\mathcal{P}})\left[\partial_{x_0}v-\sum_{i=1}^{3}\partial_{x_i}ve_i\right]&=&
\Pi_{e_4} \left[\partial_{x_0}v-\sum_{i=1}^{3}\partial_{x_i}ve_i \right]\nonumber\\&+&\sum_{i=1}^3(\Pi_{e_ie_4}+\Pi_{{\mathcal{P}}[{e_{i+1}e_{i-1}e_4}]})\left[ \partial_{x_0}v-\sum_{j=1}^{3}\partial_{x_j}ve_j \right]\nonumber\\&=& \left(\partial_{x_0}v^0-\sum_{j=1}^{3}\partial_{x_i}v^i \right)e_4\\&+&\sum_{i=1}^3 \left(\partial_{x_0}v^i+\partial_{x_i}v^0-\partial_{x_{i+1}}v^{i-1}+\partial_{x_{i-1}} v^{i+1} \right) e_ie_4.\nonumber\end{eqnarray}
 
We can associate to the operator $(\Pi_{4}+{\mathcal{P}})$ 
 the following symbol:
\begin{equation}
\sigma(\xi)=\left(\begin{array}{cccc}
\xi_0 & -\xi_1&-\xi_2&-\xi_3\\[3mm]
\xi_1&\xi_0&\xi_3&-\xi_2\\[3mm]
\xi_2&-\xi_3&\xi_0&\xi_1\\[3mm]
\xi_3&\xi_2&-\xi_1&\xi_0\end{array}\right)\end{equation}
It can be easily seen that the columns of the symbol form an orthogonal system. Therefore, we know $\operatorname{det}(\sigma(\xi))=\pm(\sum_{i=1}^4\xi_i^2)^2$ due to the multilinearity of the determinant coupled with the determinant of real orthogonal matrices being either $1$ or $-1$. This implies that the symbol is invertible for all $\xi \neq 0$ and as a result, the differential operator is elliptic by definition. Due to the connectedness of $\R^4 \setminus \{ 0 \}$ and the continuity of the determinant, we may even conclude that the sign of the determinant has to be constant and by noticing $\operatorname{det}(\sigma((1,0,0,0)))=1$, we deduce $\operatorname{det}(\sigma(\xi))=(\sum_{i=1}^4\xi_i^2)^2$ for all $\xi$. Combining the ellipticity of the Laplacian with the ellipticity of $\sigma(\xi)$, we deduce that
${\mathcal{L}}_{\frak{I}}({\frak{u}})$ is elliptic as well. This concludes the proof of Proposition \ref{invert}.~~\hfill$\Box$\par
\medskip
 We may now prove the following result:
\begin{Lm}\label{4D:la:invertible}
For any $\Theta > 0$, there exists  $\eps > 0$ so that the following holds for any $\omega_0,\frak{g}_0 $  and $\frak{q}_0$ satisfying \eqref{4D:decomp}, \eqref{4D:eq:gauge:2est}, \eqref{4D:eq:gauge:pest},  \eqref{4D:eq:gauge:ppest}:
\par
For any $\om\in ( W^{-1,p}\cap  W^{-1,p'})({\R^4},E_6)$ and 
 $\frak{g}\in   (L^p\cap L^{p^{\prime}})({\R^4},E_4)$, 
 there exists a unique ${\frak{u}}\in   \dot W^{1,p} \cap \dot W^{1,p^{\prime}} \cap L^{4p/3p-4}(\R^4,\frak{spin}(5))$, such that 
\[
(\omega,\frak{g}) = \mathcal{L}_{\frak{q}_0}({\frak{u}}) 
\]
and for some constant $C = C(\omega_0,\frak{g}_0,\Theta) > 0$, it holds
\begin{eqnarray}\label{esteta}
 \|\nabla{\frak{u}}\|_{L^{p}(\R^4)} +\|\nabla{\frak{u}}\|_{L^{p^{\prime}}(\R^4)} &\lesssim&
  \|\omega\|_{W^{-1,p}(\R^4)} +\|\omega\|_{W^{-1,p'}(\R^4)}\\&&+\|\frak{g}\|_{L^{p}(\R^4)} +\| \frak{g}\|_{L^{p^{\prime}}(\R^4)}.\nonumber 
\end{eqnarray}
\end{Lm}
\noindent
{\bf Proof of Lemma \ref{4D:la:invertible}.}\par
\noindent{\bf Claim 1.}  ${\mathcal {L}}_{\frak{I}}({\frak{u}})$ is invertible as a map between the space of functions $\frak{u}\in \dot {W}^{1,p^{\prime}} \cap L^{4p/3p-4}(\R^4, \frak{spin}(5))$  and  the space $W^{-1,{p'}}(\R^4,E_6)\times  L^{p^{\prime}} (\R^4,E_4)$\\
\noindent
 {\bf Proof of the Claim 1.}
We have seen that ${\mathcal {L}}_{\frak{I}}({\frak{u}})$ is elliptic and therefore a Cald\'eron-Zygmund operator. More precisely, let $\Gamma_4$ denote the fundamental solution of $\Delta$ on $\R^4$. Using the decomposition $\frak{u} = w + v$ as before, we have:
$$\Delta w = \omega \Longrightarrow w = \Gamma_4 \ast \omega.$$
Similarily, we write $v = v^0 e_4 + v^1 e_1e_4 + v^2 e_2e_4 + v^3 e_3e_4$ and up to replacing $e_4, e_1e_4, e_2e_4$ and $e_3e_4$ by the quaternionic basis $1, i, j$ and $k$ respectively, we see:
$$(\Pi_4 + \mathcal{P})[\partial_R v] = \frak{g} \Longleftrightarrow D^{RF}_{R} v = \frak{g},$$
where $D^{RF}_{R} = \partial_{x_0} + \partial_{x_1} \cdot i + \partial_{x_2} \cdot j + \partial_{x_3} \cdot k$ is the quaternionic Riemann-Fueter operator in $4$D. Observe that this emergence crucially limits the dimension of the domains to which this very argument could be applied. A simple calculation as outlined in the Appendix enables us to see:
$$\overline{D}^{RF}_{R} D^{RF}_{R} = \Delta,$$
where $\overline{D}^{RF}_{R} = \partial_{x_0} - \partial_{x_1} \cdot i - \partial_{x_2} \cdot j - \partial_{x_3} \cdot k$ is the conjugate operator. Therefore, we have:
$$\Delta v = \overline{D}^{RF}_{R} \frak{g}.$$
As a result, we deduce:
$$v = \Gamma_4 \ast \overline{D}^{RF}_{R} \frak{g} = \overline{D}^{RF}_{R} \frak{g} \ast \Gamma_4 = \overline{D}^{RF}_{R} \left( \frak{g} \ast \Gamma_4 \right) = \frak{g} \ast \overline{D}^{RF}_{R} \Gamma_4.$$
We highlight that the change of order in the convolution is made to emphasise explicitly the non-commutativity of elements in the Clifford algebra. Using standard Cald\'eron-Zygmund estimates for the Laplacian, we obtain:
\begin{align*}
	\| \nabla w \|_{L^{p}} 	&\lesssim \| \omega \|_{W^{-1,p^\prime}}\\
	\| \nabla v \|_{L^{p}}	&\lesssim \| \frak{g} \|_{L^{p^\prime}}.
\end{align*}
Consequently, given $\omega\in  W^{-1,{p^{\prime}}}(\R^4,E_6)$, ${\frak{g}}\in L^{p^{\prime}} (\R^4,E_4)$, there exists a unique $\frak{u}\in \dot W^{1,{p^\prime}} \cap L^{4p/3p-4}({\R^4},\frak{spin}(5))$ such that:
$${\mathcal {L}}_{\frak{I}}({\frak{u}})=(\omega,\frak{g})\,.$$
The elliptic estimates above yield in combination:
$$\|\nabla \frak{u}\|_{L^{p^\prime}}\lesssim\|\omega\|_{W^{-1,{p^\prime}}}+\|{\frak g}\|_{L^{p^{\prime}} }.$$ 
\noindent
{The claim is therefore proved.}~~\hfill$\Box$\\

\noindent
{\bf Estimate for ${\mathcal {L}}_{\frak{q}_0}({\frak{u}})-{\mathcal {L}}_{\frak{I}}({\frak{u}})$}\par
\noindent 
To generalize the invertibility to arbitrary $\frak{q}_0$, let us consider $\mathcal{L}_{\frak{q}_0}$ as a perturbation of $\mathcal{L}_{\frak{I}}$. Invertibility is ensured, if the operators are close enough by the usual perturbation-type argument. Thus, it suffices to estimate using H\"older's inequality, boundedness/compactness of the Spin-groups, Sobolev-embeddings and the $L^{4}$-estimate for $\nabla \frak{q}$:
 \begin{align*}
 	\|\frak{q}_0^{-1}(\partial_{x_i}\frak{q}_0)\frak u-\frak{u}\frak{q}^{-1}_0\partial_{x_i}\frak{q}_0\|_{L^{p^\prime}}	&\lesssim \| \frak{q}_0^{-1} \|_{L^\infty} \|\nabla \frak{q}_0\|_{L^4}\|\frak{u}\|_{L^{4p/3p-4}} \\
																						&\lesssim \|\nabla \frak{q}_0\|_{L^4}\|\nabla \frak{u}\|_{L^{p^\prime}} \\
																						&\lesssim \Theta \varepsilon \cdot \|\nabla \frak{u}\|_{L^{p^\prime}}.
 \end{align*}
 Using this inequality, we conclude:
 \begin{align*}
 	\|\partial_{x_i}\left(\frak{q}_0^{-1}(\partial_{x_i}\frak{q}_0)\frak {u}-\frak{u}\frak{q}^{-1}_0\partial_{x_i}\frak{q}_0\right)\|_{W^{-1,p^\prime}}	&\leq \|\frak{q}_0^{-1}(\partial_{x_i}\frak{q}_0)\frak u-\frak{u}\frak{q}^{-1}_0\partial_{x_i}\frak{q}_0\|_{L^{p^\prime}} \\
																											&\lesssim \Theta \varepsilon \cdot \|\nabla \frak{u}\|_{L^{p^\prime}}.
 \end{align*}
Choosing $\varepsilon>0$ small enough (depending on $\Theta$), we obtain that ${\mathcal {L}}_{\frak{q}_0}({\frak{u}})$ is an invertible map from $\dot{W}^{1,p'}(\R^4, \frak{spin}(5))$ to $W^{-1,{p^{\prime}}}(\R^4,E_6)\times  L^{p^{\prime}} (\R^4,E_4).$ 

\medskip
\noindent{\bf Claim 2.} Assuming now $\omega \in  ( W^{-1,p}\cap  W^{-1,p'})({\R^4},E_6)$ as well as $\frak{g}\in  (L^p\cap L^{p^{\prime}})({\R^4},E_4)$, we show that
the unique solution  ${\frak{u}}$ of $(\omega,\frak{g}) = \mathcal{L}_{\frak{q}_0}({\frak{u}})$ lies in $\dot W^{1,p}(\R^4)$. \par
\noindent
{\bf Proof of Claim 2:} Firstly, due to $\nabla \frak{u} \in L^{p^\prime}$, we know that we may choose $\frak{u}$ by Sobolev-embeddings and the density of Schwartz functions in the following way:
$$\frak{u} \in L^{4p/3p-4}(\R^4)$$
We have been using this observation implicitly before. As previously, in order to bootstrap, it suffices to deduce improved integrability of $\frak{q}_0^{-1} \partial_{x_l} \frak{q}_0 \frak{u}$, as this implies improved integrability of $\nabla \frak{u}$ by means of elliptic estimates. The same estimates immediately apply to $\frak{u} \frak{q}_0^{-1} \partial_{x_l} \frak{q}_{0}$, meaning that there is no issue in merely establishing estimates for $\frak{q}_0^{-1} \partial_{x_l} \frak{q}_0 \frak{u}$ for brevity's sake. By the considerations in \eqref{contembedlppprime}, it suffices to show that $\nabla \frak{u} \in L^{q}$ for some $q > 4$, because then, by interpolation, $\nabla \frak{u} \in L^{p^\prime} \cap L^{q}$ and we could thus conclude that $\frak{u} \in L^{\infty}$ completely analogous to \eqref{contembedlppprime} leading to $\frak{q}_0^{-1} \partial_{x_l} \frak{q}_0 \frak{u} \in L^{p}$, which immediately establishes $\nabla \frak{u} \in L^{p}$. Therefore, Claim 2 would be proven in the process. \par
We argue by a bootstrap argument: Assume that $\frak{u} \in L^{r}$ for some $4 > r \geq \frac{4p}{3p-4}$. In this case, H\"older's inequality implies:
$$\| \frak{q}_0^{-1} \partial_{x_l} \frak{q}_0 \frak{u} \|_{L^{t}} \lesssim \| \nabla \frak{q}_{0} \|_{L^{p}} \| \frak{u} \|_{L^{r}},$$
for $\frac{1}{t} = \frac{1}{p} + \frac{1}{r} > \frac{1}{p} + \frac{1}{4} > \frac{1}{4}$. Observe that $\frac{4}{3} \leq t < 4 < p$ by the inequalities satisfied by $r$. We conclude due to the elliptic estimates as in Claim 1 and the identity $\mathcal{L}_{\frak{q}_{0}}(\frak{u}) = (\omega, \frak{g})$:
$$\nabla \frak{u} \in L^{t}.$$
This implies by Sobolev-embeddings that $\frak{u} \in L^{4t/4-t}$. Thus, if we define $\tilde{r} = \frac{4t}{4 - t}$, we observe:
$$\frac{1}{\tilde{r}} = \frac{1}{t} - \frac{1}{4} = \frac{1}{r} + \frac{1}{p} - \frac{1}{4},$$
which implies that the reciprocal values are decreasing by a constant amount with each iterating step, due to $p > 4$. Therefore, after finitely many steps (the number of which depends only on $p$), we have:
$$\frac{1}{\tilde{r}} < \frac{1}{4}\ \Rightarrow\ \tilde{r} > 4$$
This implies, by the previously outlined argument, that $\nabla \frak{u} \in L^{p}(\R^4, \frak{spin}(5))$, finishing the proof of Claim 2. Observe that by keeping track of the estimates, we may deduce the $L^p$-part of the inequality \eqref{esteta}. Therefore, the Lemma is proven.
\qed
\par
\medskip
\par
 
\noindent{\bf Proof of   Proposition \ref{4D:pr:VepseqUeps} continued.}\\
  \noindent
 For $\eps = \eps(\Theta) > 0$ chosen small enough and for any $(\omega_0,\frak{g}_0) \in \mathcal{V}_{\eps,\Theta}$, the local inversion theorem applied to ${\mathcal{N}}_{\frak{q}_0}$ implies the existence of some $\delta > 0$ (that might depend on $(\omega_0,\frak{g}_0)$) such that, for every $(\omega,\frak{g}) \in {\mathcal{U}}_\eps$ with

\begin{eqnarray}
&& \|\omega - \omega_0 \|_{W^{-1,p}(\R^4)}+\|\omega - \omega_0 \|_{W^{-1,p^{\prime}}(\R^4)} < \delta\label{4D:estgauge5}\\[5mm]
&& \|\frak{g} - \frak{g}_0 \|_{L^{p}(\R^4)}+\|\frak{g}- \frak{g}_0 \|_{L^{p^{\prime}}(\R^4)} < \delta, \label{4D:estgauge6}
\end{eqnarray}
we surely find $\frak{q} = \frak{q}_0 e^{\frak{u}} \in \dot{W}^{1,p} \cap \dot{W}^{1,p^{\prime}}(\R^4 )$, such that $ \frak{q}-\frak{I}  \in  L^{4p/3p-4}(\R^4)$ and \eqref{4D:decomp} is satisfied.

It remains to prove the estimates \eqref{4D:eq:gauge:2est}, \eqref{4D:eq:gauge:pest} and \eqref{4D:eq:gauge:ppest}. They will be an immediate consequence of the following lemma, together with sufficiently small chosen $\varepsilon, \delta > 0$:
\begin{Lm}\label{4D:la:qgaugeuniformest}
There exist $\Theta > 0$ and $\sigma > 0$, such that, whenever $\frak{q} \in \dot{W}^{1,p} \cap \dot{W}^{1,p^{\prime}}(\R^4)$ with $ \frak{q}-\frak{I}  \in  L^{4p/3p-4}(\R^4)$ satisfying \eqref{4D:decomp} is given, and it holds:
\begin{equation}
  \|\nabla\frak{q}\|_{L^4(\R^4)}\leq \sigma,
\end{equation}
then  the estimates in \eqref{4D:eq:gauge:2est}, \eqref{4D:eq:gauge:pest} and \eqref{4D:eq:gauge:ppest} hold true as well.
\end{Lm}
\noindent
{\bf Proof of Lemma \ref{4D:la:qgaugeuniformest}.} 
Let $(\omega,\frak{g}) \in {\mathcal{U}}_\eps$ satisfy \eqref{4D:estgauge5} and  \eqref{4D:estgauge6} and let    $\frak{q} = \frak{q}_0 e^{\frak{u}} \in \dot{W}^{1,p} \cap \dot{W}^{1,p^{\prime}}(\R^4 )$, such that  $ \frak{q} -\frak{I} \in  L^{4p/3p-4}(\R^4)$  and \eqref{4D:decomp} is satisfied. We first consider the following Hodge decomposition of $\frak{q}^{-1}d\frak{q}$:
\begin{equation}\label{4D:hodge3}
\frak{q}^{-1}d\frak{q}=d\Gamma_{\frak{q}}+d^*Y_{\frak{q}}
\end{equation}
where $Y_{\frak{q}}\in \Omega^2(\R^4 )$, $Y_{\frak{q}}=\sum_{0\le i<j\le 3}Y_{\frak{q}}^{ij} dx_i\wedge dx_j$ is a differential $2$-form and $\Gamma_{\frak{q}}$ a $0$-form, i.e. a function. We denote $d^*Y_{\frak{q}}=\sum_{i=0}^{3}y^i_{\frak{q}}dx_i$\footnote{We recall that  $d^{*}=(-1)^{n(k-1)+1}\ast \ d\ \ast$, $\ast$ is the Hodge operator. If $\xi=\sum_{0\le i<j\le 3}\xi_{ij} dx_i\wedge dx_j$ then
$d^*\xi=-(\alpha_0 dx_0+\alpha_1 dx_1+\alpha_2 dx_2+\alpha_3 dx_3)$ where:
\begin{eqnarray*}
\alpha_0&=& \partial_{x_1}\xi_{01}+\partial_{x_2}\xi_{02}+\partial_{x_3}\xi_{03}\\
\alpha_1&=& -\partial_{x_0}\xi_{01}+\partial_{x_2}\xi_{12}+\partial_{x_3}\xi_{13}\\
\alpha_2&=&- \partial_{x_0}\xi_{02}-\partial_{x_1}\xi_{12}+\partial_{x_3}\xi_{23}\\
\alpha_3&=& -\partial_{x_2}\xi_{23}-\partial_{x_1}\xi_{13}-\partial_{x_0}\xi_{03}
\end{eqnarray*}
}
for brevity. We may choose $\Gamma_{\frak{q}}$ and $Y_{\frak{q}}$ as follows:
\begin{eqnarray}\label{4D:YG2}
\Gamma_{\frak{q}}&=&(-\Delta)^{-1}d^*(\frak{q}^{-1}d\frak{q})\label{Gamma2}\\
Y_{\frak{q}}&=&(-\Delta)^{-1}d(\frak{q}^{-1}d\frak{q})\label{Y2}
\end{eqnarray}
In particular, we then have $dY_{\frak{q}}=0$ and $d^*\Gamma_{\frak{q}}=0$, i.e. exactness and coexactness respectively.\par
Due to \eqref{4D:YG2}, it follows that:
\begin{equation}\label{estGamma6}
(-\Delta) \Pi_6(\Gamma_{\frak{q}}) = \Pi_6(-\Delta \Gamma_{\frak{q}})=\Pi_6(d^*(\frak{q}^{-1}d\frak{q}))=-\omega.
\end{equation}
Therefore for every $r\in [p',p]$ we have:
\begin{equation}\label{4D:estgradientgamma}
\|\Pi_6(\nabla \Gamma_{\frak{q}})\|_{L^r}\lesssim \|\omega\|_{W^{-1,r}}
\end{equation}
 
From \eqref{4D:hodge3}, it follows that 
\begin{eqnarray}\label{4D:regY3}
-\Delta Y_{\frak{q}}&=&d(\frak{q}^{-1}d\frak{q})=d\frak{q}^{-1}\wedge d\frak{q}.
\end{eqnarray}

Using \eqref{4D:regY3}, it follows that  $\nabla Y_q\in L^{r}(\R^4)$ and due to the compensation result in Lemma \ref{prel:compensation}:
\begin{equation}\label{estgradientY}
\|\nabla Y_{\frak{q}}\|_{L^r}\lesssim \| d\frak{q}\|_{L^{4}(\R^4)}\| d\frak{q}\|_{L^{r}(\R^4)}\le \sigma \| d\frak{q}\|_{L^{r}(\R^4)}.
\end{equation}

%Now we observe  that $\frak{q}^{-1}d\frak{q}$ is {\em purely  imaginary}, namely it is a linear combination of elements in
%$\frak{spin}(5)$, since $|\frak{q}|=1$. Hence:
%\begin{equation}\label{pi4P}
% (\Pi_{4}+{\mathcal{P}})(\sum_{i=1}^3\frak{q}^{-1}\partial_{x_i}\frak{q}e_i)=(\Pi_4+{\mathcal{P}})(\sum_{i=1}^3\frak{q}^{-1}\partial_{x_i}\frak{q})e_i.\end{equation}
 By inserting  \eqref{4D:hodge3}, it follows that (we write $D = \partial_{R}$ for the moment for brevity's sake):
 \begin{eqnarray}\label{pi4PD}
 \frak{g}=(\Pi_4+{\mathcal{P}}){\mathcal{D}}(\frak{q})&=&(\Pi_{4}+{\mathcal{P}})(D\Gamma_{\frak{q}})\nonumber\\
 &+&(\Pi_{4}+{\mathcal{P}})(y^0_{\frak{q}}-\sum_{i=1}^3y^i_{\frak{q}}e_i).
 \end{eqnarray}
 Therefore:
 \begin{equation}
 (\Pi_{4}+{\mathcal{P}})(D\Gamma_{\frak{q}})=\frak{g}-(\Pi_{4}+{\mathcal{P}})(y^0_{\frak{q}}-\sum_{i=1}^3y^i_{\frak{q}}e_i)
 \end{equation}
 Observe that $d\Gamma_{\frak{q}}\in \frak{spin}(5)$, since $\frak{q}^{-1}d\frak{q}\in \frak{spin}(5)$. Therefore:
$$d\Gamma_{\frak{q}}=(\Pi_4+\Pi_6)(d\Gamma_{\frak{q}}).$$
Hence:
\begin{eqnarray}
 (\Pi_{4}+{\mathcal{P}})(D\Pi_4(d\Gamma_{\frak{q}}))&=&d\frak{g}-d(\Pi_{4}+{\mathcal{P}})(y^0_{\frak{q}}-\sum_{i=1}^3y^i_{\frak{q}}e_i)\\
 &-&(\Pi_{4}+{\mathcal{P}})(D\Pi_6(d\Gamma_{\frak{q}}))
 \end{eqnarray}
  Since the operator $(\Pi_{4}+{\mathcal{P}})\circ D$ is invertible by the arguments in Claim 1 of the proof of Lemma \ref{4D:la:invertible} above, we find:
  \begin{eqnarray}\label{4D:P4dGamma}
 \Pi_4(d\Gamma_{\frak{q}})&=&((\Pi_4+{\mathcal{P}})\circ D)^{-1}d\frak{g} \nonumber \\
 &+&((\Pi_4+{\mathcal{P}})\circ D)^{-1}d(\Pi_{4}+{\mathcal{P}})\left( y^0_{\frak{q}}-\sum_{i=1}^3y^i_{\frak{q}}e_i \right)\nonumber \\
 &+&((\Pi_4+{\mathcal{P}})\circ D)^{-1}\left[ (\Pi_{4}+{\mathcal{P}})(D\Pi_6(d\Gamma_{\frak{q}}))\right].
\end{eqnarray}
By using \eqref{4D:P4dGamma}, we get:

 \begin{eqnarray}\label{4D:estgradientgamma2}
 \|\Pi_{4} (d\Gamma_{\frak{q}})\|_{L^r}&\lesssim& \|\frak{g}\|_{L^r}+\|d^*Y_{\frak{q}}\|_{L^r}+\|\omega\|_{W^{-1,r}}\nonumber
\\&\lesssim&  \|\frak{g}\|_{L^r}+\sigma \| d\frak{q}\|_{L^{r}}+\|\omega\|_{W^{-1,r}}
\end{eqnarray}
 Combining \eqref{4D:hodge3},  \eqref{4D:estgradientgamma} and \eqref{4D:estgradientgamma2}, we get the following estimate:
 \begin{eqnarray}\label{estgradientq}
 \|d\frak{q}\|_{L^r}&\lesssim& \|d\Gamma_{\frak{q}}\|_{L^r}+\|d^*Y_{\frak{q}}\|_{L^r}\\
 &\lesssim&\|\Pi_4(d\Gamma_{\frak{q}})\|_{L^r}+\|\Pi_6(d\Gamma_{\frak{q}})\|_{L^r}+\|d^*Y_{\frak{q}}\|_{L^r}\nonumber\\
 &\le & C(\|\frak{g}\|_{L^r}+\sigma \| d\frak{q}\|_{L^{r}}+2\|\omega\|_{W^{-1,r}}+ \sigma \| d\frak{q}\|_{L^{r}}).
 \end{eqnarray}  
Choosing $\Theta :=\frac{C}{1-2C\sigma}$, we finally arrive at the desired inequality: 
$$
 \|d\frak{q}\|_{L^r}\le \Theta(\|\omega\|_{W^{-1,r}}+ \| \frak{g}\|_{L^r}).$$
 This concludes the proof of lemma~\ref{4D:la:qgaugeuniformest}.~~\hfill$\Box$\\
 
 \noindent
 {\bf End of the proof of Proposition \ref{4D:pr:VepseqUeps}}
 
\noindent
Thanks to Lemma~\ref{4D:la:qgaugeuniformest}, the openness property (iv.) is proven. Proposition~\ref{4D:pr:VepseqUeps} is thus established. \qed
\par
\medskip

\subsection{Improved Integrability}
 
We are now going to finish the proof of Theorem \ref{mainth4D}. Before we start, however, let us briefly recall the definition of the gauge operator and the conditions: Let $\frak{f}\in L^{4/3}(\R^4)$ be a solution of
\begin{eqnarray}\label{maineq3}
\partial_{x_0}(\frak{q}\frak{f}) -\sum_{i=1}^{3} \partial_{x_i}(\frak{q}e_i \frak{f} )&=&
\frak{q}\left(\beta e_4+{\mathcal{D}}(\frak{q})\right) \frak{f}.
\end{eqnarray}
If $\|\beta\|_{L^{(4,2)}(\R^4)}\le \varepsilon$ (this is the required regularity assumption for our arguments, the corresponding $L^{4}$-estimate follows immediately) for some $\varepsilon > 0$ sufficiently small, then there exists $\frak{q}\in \dot W^{1,4}(\R^4)$   such that
\begin{equation}\label{bt1}
{\mathcal{N}}(\frak{q})=(0,-\beta e_4)
\end{equation}
with
\begin{equation}\label{bt2}
\|\nabla\frak{q}\|_{L^4(\R^4)}\le \Theta \|\beta\|_{L^4(\R^4)}
\end{equation}
This is what we have proven in the last subsection. Here, $\mathcal{N}$ denotes the following gauge operator:
 \begin{eqnarray}\label{est9}
&& {\mathcal{N}} \colon \dot W^{1,4}(\R^4 ,Spin(5))\to  W^{-1,4}(\R^4 ,E_6)\times  L^4(\R^4,E_4) \\[5mm]
&&~~~~~ \frak{q}\mapsto\left(\Pi_6\left(\sum_{i=0}^3(\partial_{x_i}(\frak{q} ^{-1}\partial_{x_i} \frak{q})\right), ( \Pi_{4}+{\mathcal{P}})({\mathcal{D}}(\frak{q}))\right) \nonumber
\end{eqnarray}
In order to avoid worrying about signs, we shall from now on work with $\beta e_4$ instead of $-\beta e_4$. This can be achieved by replacing $\beta$ by $-\beta$ and does not affect the argument in any meaningful way.\\
In particular, it follows from \eqref{bt1} that:
\begin{equation}\label{bt3}
(\Pi_{4}+{\mathcal{P}}){\mathcal{D}}(\frak{q})=\beta e_4.
\end{equation}
Namely, if $\beta = (\beta^0, \beta^1, \beta^2, \beta^3)$:
\begin{eqnarray}\label{bt3bis}
&&\Pi_{e_4}({\mathcal{D}}(\frak{q}))=\beta^0\\
&&(\Pi_{e_ie_4}+\Pi_{{\mathcal{P}}({e_{i+1}e_{i-1}e_4})}({\mathcal{D}}(\frak{q}))=\beta^i e_ie_4
\end{eqnarray}
An important step in the proof of our regularity result stems from the observation that the solution of this type of problem can be easily computed directly. This can be exploited to obtain further information and stronger integrability properties as seen below:
\begin{Lm}\label{4D:regP}
Under the above assumptions, we have
\begin{equation}\label{bt4}
\Pi_{6}({\mathcal{D}}(\frak{q})),~~\Pi_{e_{i+1}e_{i-1}e_4}({\mathcal{D}}(\frak{q}))\in L^{(4,1)}(\R^4).
\end{equation}
\end{Lm}
We first prove a related result concerning the linearized operator ${\mathcal{L}}_{\frak{I}}$. For convenience's sake, given $\frak{u} = w + v \in E_6 \oplus E_4$, we set:
\begin{equation}\label{DV}Dv := \partial_R v=\partial_{x_0}v-\sum_{i=1}^3\partial_{x_i}v e_i.\end{equation}
 
The result in Lemma \ref{4D:regP} has an infinitesimal analogue for the differential which is in fact the key element required to prove it:

\begin{Lm}\label{4D:regP2}
Let $\frak{u}=w+v\in E_6 \oplus E_4$ be such that
\begin{equation}\label{bt5}
{\mathcal{L}}_{\frak{I}}(\frak{u})=(\Delta w,\Pi_{e_{4}}(Dv),(\Pi_{e_ie_4}+\Pi_{{\mathcal{P}}({e_{i+1}e_{i-1}e_4})}(Dv))=(0,\beta e_4)
\end{equation}
Then for all $i=1,2,3$ we have
\begin{equation}\label{bt6}
\Pi_{e_{i+1}e_{i-1}e_4}(Dv)=0,
\end{equation}
and therefore:
$$\Pi_{{\mathcal{P}}({e_{i+1}e_{i-1}e_4})}(Dv)=0.$$
\end{Lm}

The key idea behind the proof is the use of explicit representations of the solution $\frak{u}$.\\

\noindent
{\bf Proof of Lemma \ref{4D:regP2}.}\\

We write $v=v^0e_4+v^1e_1e_4+v^2e_2e^4+v^{3}e_3e_4$ as 
$v=(v^0,v^{\prime})$ where $v^{\prime}=(v^1,v^2,v^3)$ and similarily $x=(x_0,x^\prime)$, where $x^\prime=(x_1,x_2,x_3)$. We observe that $Dv$ can be computed as follows:
\begin{eqnarray}\label{bt7}
Dv&=&(\partial_{x_0}v^0-\div_{x^{\prime}}v^{\prime})e_4+\sum_{i=1}^3(\partial_{x_i}v^{0}+\partial_{x_0}v^{i})e_ie_4\\
&+&(\partial_{x_3}v^{2}-\partial_{x_{2}}v^{3})e_{2}e_{3}e_{4}+(\partial_{x_1}v^{3}-\partial_{x_{3}}v^{1})e_{3}e_{1}e_{4}+(\partial_{x_2}v^{1}-\partial_{x_{1}}v^{2})e_{1}e_{2}e_{4}\nonumber
\end{eqnarray}
Therefore, we may express $Dv$ in the following form:
\begin{equation}\label{bt8}
Dv=\left(\begin{array}{c}
\partial_{x_0}v^0-\div_{x^{\prime}}v^{\prime}\\
\nabla_{x^{\prime}}v^{0}+\partial_{x_0}v^{\prime}\\
-{\rm curl}_{x^{\prime}} v^{\prime}\end{array}\right)
\end{equation}
\noindent
We want to find the solution $v\in \dot W^{1,4}(\R^4)$ of the following system of equations: 
\be
\label{bt8-aa}Dv=\left(\begin{array}{c}
\beta^0\\\beta^1\\ \beta^2\\ \beta^3\\0\\0\\0\end{array}\right)\ .
\ee
$\mathbf 1)  $ Assume that $\beta\in {\mathcal{S}}(\R^4)$. We show the existence of a smooth solution $v$ and look for a-priori estimates.\\

\noindent
First of all, we notice that:
\begin{align}
	\Delta v^0	&= \partial_{x_0} \left( \partial_{x_0} v^{0} \right) + \partial_{x_1} \left( \partial_{x_1} v^{0} \right) + \partial_{x_2} \left( \partial_{x_2} v^{0} \right) + \partial_{x_3} \left( \partial_{x_3} v^{0} \right) \notag \\
			&= \partial_{x_0} \left( \partial_{x_0} v^{0} \right) + \div_{x'}  \nabla_{x'} v^0 \notag \\
			&= \partial_{x_0} \left( \beta^0 + \div_{x'} v' \right) + \div_{x'} \left( \beta' - \partial_{x_0} v' \right) \notag \\
			&= \div \beta,
\end{align}
and thus:
$$v^0(x) := (-\Delta)^{-1} \left( \div \beta \right)(x) = -\int_{\R^4}\div(\beta)(y)|x-y|^{-2} dy$$
Our goal is now to arrive at similar expressions for $v^{j}$ for all $j = 1,2,3$. To achieve this, we observe that for any such $j$:
\begin{align}
	\partial_{x_{0}} \beta^{j} - \partial_{x_{j}} \beta^{0}	&= \partial_{x_{0}} \left( \partial_{x_{j}} v^0 + \partial_{x_0} v^{j} \right) - \partial_{x_{j}} \left( \partial_{x_0} v^0 - \div_{x'} v' \right) \notag \\
									&= \partial_{x_{0}}^{2} v^{j}   + \sum_{k\neq j} \partial_{x_{j}} \partial_{x_{k}} v^{k} \notag \\
									&= \partial_{x_{0}}^{2} v^{j}   + \sum_{k\neq j} \partial_{x_{k}} \partial_{x_{j}} v^{k} \notag \\
									&= \partial_{x_{0}}^{2} v^{j}   + \sum_{k\neq j} \partial_{x_{k}} \partial_{x_{k}} v^{j} \notag \\
									&= \Delta v^{j},
\end{align}
where we used that $\partial_{x_{j}} v^{k} = \partial_{x_{k}} v^{j}$ for all $j \neq k$ should hold by the third set of equations in \eqref{bt8} (namely ${\rm curl}_{x^{\prime}} v^{\prime}=0$). Thus, we also know:
$$v^{j}(x) := (-\Delta)^{-1}\left( \partial_{x_{j}} \beta^0 - \partial_{x_0} \beta^{j} \right)(x) = -\int_{\R^4}\left( \partial_{x_{j}} \beta^0(y) - \partial_{x_0} \beta^{j}(y) \right) |x-y|^{-2} dy$$
We observe that $v$ obtained this way clearly satisfies the desired $L^4$-estimate by the usual Calderon-Zygmund inequality. Consequently, we merely have to verify that this solution does indeed solve the equation \eqref{bt8-aa}. Since this is done by direct computations, let us only present the computations in the case of the second set of equations in \eqref{bt8}:
\begin{align}
	\partial_{x_{j}} v^{0} + \partial_{x_0} v^{j}	&= -\partial_{x_{j}} \left( \Gamma \ast \div \beta \right) + \partial_{x_{0}} \left( \Gamma \ast \left( \partial_{x_j} \beta^0 - \partial_{x_0} \beta^j \right) \right) \notag \\
								&= \Gamma \ast \left( -\partial_{x_j} \partial_{x_0} \beta 0 - \partial_{x_{j}}^2 \beta^j - \sum_{k \neq j} \partial_{x_{j}} \partial_{x_k} \beta^k + \partial_{x_{0}} \partial_{x_{j}} \beta^{0} - \partial_{x_{j}}^2 \beta^j \right) \notag \\
								&= \Gamma \ast \left( -\Delta \beta^{j} - \sum_{k \neq j} \partial_{x_{k}} \left( \partial_{x_{j}} \beta^{k}- \partial_{x_{k}} \beta^{j} \right) \right) \notag \\
								&= \Gamma \ast \left( - \Delta \beta^{j} \right) = \beta^{j},
\end{align}
where we denote by $\Gamma$ the fundamental solution of the Laplacian $-\Delta$ in $4$D and we used $\mathbf{curl}_{x'} \beta' = 0$. This computation is valid for any $j = 1,2,3$. This shows that the second set of equations in \eqref{bt8} holds true and the other two sets of equations may be checked completely analogously and are omitted here.\\

\noindent
$\mathbf{2) }$ The general case, i.e. the case of $\beta \in L^{4}(\R^4;\operatorname{span}_{\R} \{ e_0, e_1, e_2, e_3 \})$ satisfying the vanishing curl condition, can be dealt with by approximation. Notice that any such $\beta$ can be approximated by Schwartz functions or smooth, compactly supported functions for which the previous computations hold. Then, the uniformity of the estimates on the gradient of $v$ leads to the desired conclusion.\\

  A particular special case is when $\beta = \partial_{L} \alpha$ for some $\alpha \in \dot{W}^{1,(4,2)}(\R^4)$ real-valued. Keep in mind that:
 \begin{equation}
 	\partial_{L} \alpha = \partial_{x_0} \alpha - \partial_{x_{1}} \alpha \cdot e_1 - \partial_{x_{2}} \alpha \cdot e_2 - \partial_{x_{3}} \alpha \cdot e_3.
 \end{equation}
In fact, in this case, we may find an even more explicit representation of the solution $v$. Indeed, by the vanishing curl assumption on $v'$, it is natural to look for solutions:
  \begin{equation}
  \label{naturalpartialgradient}
  	v' = \nabla_{x'} \varphi
  \end{equation}
  Inserting this expression into the second set of equations in \eqref{bt8}, we find:
  $$\nabla_{x'} v^{0} + \partial_{x_{0}} \nabla_{x'} \varphi = -\nabla_{x'} \alpha,$$
  where we remember that we currently assume $\beta = \partial_{L} \alpha$. This immediately yields:
  $$\nabla_{x'} \left( v^{0} + \partial_{x_{0}} \varphi + \alpha \right) = 0,$$
  which would be satisfied, if for instance:
  \begin{equation}
  \label{helpfulassumption}
  	v^{0} + \partial_{x_{0}} \varphi + \alpha = 0.
  \end{equation}
  It remains to check whether the first equation in \eqref{bt8} can hold true. Inserting yields:
  $$\partial_{x_{0}} v^{0} - \div_{x'} v' = \partial_{x_{0}} \alpha,$$
  which, by using the identity from \eqref{helpfulassumption} in the following form:
  $$v^{0} = -\alpha - \partial_{x_{0}} \varphi,$$
  further reduces to:
  \begin{equation}
  	- \partial_{x_{0}}^{2} \varphi - \div_{x'} \nabla_{x'} \varphi = - \Delta \varphi = 2\partial_{x_{0}} \alpha.
  \end{equation}
  Therefore:
  $$\varphi := 2 (-\Delta)^{-1} \partial_{x_{0}} \alpha,$$
  and $v^{0}, v'$ can now be computed from \eqref{naturalpartialgradient} and \eqref{helpfulassumption}. The desired estimates are evident from our computations and using that $\nabla \alpha \in L^{(4,2)}$, i.e. the gradient of $\alpha$ possesses $L^{4}$-integrability. Notice that the formula provides the same result as in the previous computation for general $\beta$.
~~\hfill$\Box$\\

It should be noted that the arguments in the previous section do not make use of the dimension of the domain being $4$ in any meaningful way, besides ensuring that a connection to the gauge operator $\mathcal{N}$ exists. Indeed, the very same arguments could be applied in other dimensions, in particular the construction of a curl-free solution of a system of PDEs.\\

\noindent
{\bf Proof of Lemma \ref{4D:regP}.}
We argue in different steps:\\

\noindent
{\bf Step 1.} We    consider the Hodge decomposition of 
\begin{equation}\label{hodge}
\frak{q}^{-1}d\frak{q}=d\Gamma_{\frak{q}}+d^*Y_{\frak{q}}
\end{equation}
where $Y_{\frak{q}}\in \Omega^2(\R^4 )$, $Y_{\frak{q}}=\sum_{0\le i<j\le 3}Y_{\frak{q}}^{ij} dx_i\wedge dx_j$ is a differential $2$-form and $\Gamma_{\frak{q}}$ is a $0$-form or function. We denote as before $d^*Y_{\frak{q}}=\sum_{i=0}^{3}y^i_{\frak{q}}dx_i$. Notice that once again, we can choose $\Gamma_{\frak{q}}$ and $Y_{\frak{q}}$ as follows:
\begin{eqnarray}\label{YG}
\Gamma_{\frak{q}}&=&(-\Delta)^{-1}d^*(\frak{q}^{-1}d\frak{q})\label{Gamma}\\
Y_{\frak{q}}&=&(-\Delta)^{-1}d(\frak{q}^{-1}d\frak{q})\label{Y}
\end{eqnarray}
In particular, we have $dY_{\frak{q}}=0$ and $d^*\Gamma_{\frak{q}}=0$, i.e. exactness and coexactness respectively. Moreover $\nabla Y_{\frak{q}},
\nabla\Gamma_{\frak{q}}\in L^4(\R^4).$

\noindent 
Due to \eqref{hodge}, it follows that 
\begin{eqnarray}\label{regY2}
-\Delta Y_{\frak{q}}&=&d(\frak{q}^{-1}d\frak{q})=d\frak{q}^{-1}\wedge d\frak{q}\nonumber\\
&=&d\frak{q}^{-1} \frak{q} \wedge \frak{q}^{-1}d\frak{q}= -(\frak{q}^{-1}d\frak{q}\wedge \frak{q}^{-1}d\frak{q})\nonumber\\
&=&-(d\Gamma_{\frak{q}}+d^*Y_{\frak{q}}
)\wedge (d\Gamma_{\frak{q}}+d^*Y_{\frak{q}})
\in L^{4}\cdot L^{4}\hookrightarrow L^{2}(\R^4)
\end{eqnarray}
From \eqref{regY2}, it follows that $\nabla^2 Y_{\frak{q}} \in L^2$ and by generalized Sobolev embeddings therefore  $\nabla Y_{\frak{q}}\in L^{(4,2)}(\R^4)$.
\par
Since $\Pi_{6}(d^{*}(\frak{q}^{-1}d\frak{q}))=0$ by the choice of $\frak{q}$ using \eqref{bt1}, we deduce from \eqref{YG} that $\Pi_{6}(\Delta \Gamma_{\frak{q}})=0$ and since $\nabla \Gamma_{\frak{q}}\in L^4(\R^4)$, this leads us to: 
\begin{equation}\label{regG2}
\Pi_{6}( \nabla\Gamma_{\frak{q}})=0
\end{equation}

\noindent   
{\bf Step 2.}  
 Next, we have, by using $D$ as in Lemma \ref{4D:regP2}:
 \begin{eqnarray}\label{estPi4PD}
\beta e_4&=& (\Pi_4+{\mathcal{P}}){\mathcal{D}}(\frak{q})=(\Pi_{4}+{\mathcal{P}})(D\Gamma_{\frak{q}})\nonumber\\
 &+&(\Pi_{4}+{\mathcal{P}})(y^0_{\frak{q}}-\sum_{i=1}^3y^i_{\frak{q}}e_i).
 \end{eqnarray}
 Since  $\frak{q}^{-1}d\frak{q}$ is {\em purely  imaginary}, namely it is a linear combination of elements in
$\frak{spin}(5)$, and $\Pi_6(d\Gamma_{\frak{q}})=0$ due to \eqref{estPi4PD}, we find: 
%$$(\Pi_{4}+{\mathcal{P}})(\sum_{i=1}^3\frak{q}^{-1}\partial_{x_i}\frak{q}e_i)=(\Pi_{4}+{\mathcal{P}})(\sum_{i=1}^3\frak{q}^{-1}\partial_{x_i}\frak{q})e_i.$$   
 \begin{eqnarray}\label{estDgamma}
 (\Pi_{4}+{\mathcal{P}})(Dd\Gamma_{\frak{q}}) &=&d\beta e_4-d(\Pi_{4}+{\mathcal{P}})(y^0_{\frak{q}}-\sum_{i=1}^3y^i_{\frak{q}}e_i) \end{eqnarray}
From \eqref{estDgamma} and the invertibility of $(\Pi_{4}+{\mathcal{P}})\circ D$, it follows that
  \begin{eqnarray}\label{dGamma}
 d\Gamma_{\frak{q}}&=& ((\Pi_{4}+{\mathcal{P}})\circ D)^{-1} (d\beta e_4)\\&-&((\Pi_{4} +{\mathcal{P}})\circ D)^{-1}d(\Pi_4+{\mathcal{P}})(y^0_{\frak{q}}-\sum_{i=1}^3y^i_{\frak{q}}e_i) .\nonumber
 \end{eqnarray}

  Now we set $\tilde Y_{\frak{q}}:= ((\Pi_{4}+{\mathcal{P}})\circ D)^{-1}(\Pi_4+{\mathcal{P}})(y^0_{\frak{q}}-\sum_{i=1}^3y^i_{\frak{q}}e_i)$ and
  let $v$ be such that $(\Pi_{4}+{\mathcal{P}})Dv=\beta  e_4$. Existence is justified by ellipticity and using the connection to the Riemann-Fueter operator introduced in the previous subsection. Observe that by elliptic estimates, we have $\nabla v\in  L^{(4,2)}$ since  $\beta\in L^{4,2}$. This is the key-point where we need that   $\beta\in L^{4,2}$.
  Therefore  $\nabla \Gamma_{\frak{q}} \in L^{(4,2)}$ as well with
  \begin{equation}\label{nablaG}
  \| \nabla \Gamma_{\frak{q}} \|_{L^{(4,2)}} \lesssim \| \nabla \frak{q} \|_{L^{4}}^{2} + \| \beta \|_{L^{(4,2)}}.\end{equation}
   We estimate:
  \begin{eqnarray}\label{estwedge}
  (d\Gamma_{\frak{q}} 
)\wedge (d\Gamma_{\frak{q}} )&=& dv\wedge dv+
 d v\wedge d\tilde Y_{\frak{q}}\nonumber\\&+&
 d\tilde Y_{\frak{q}}\wedge d v+
 d\tilde Y_{\frak{q}}\wedge d\tilde Y_{\frak{q}}
 \end{eqnarray}
Now observe that all terms are products of functions in $L^{(4,2)}$. Therefore, the product lies in $L^{2,1}$ by the Lorentz-H\"older inequality. Similarly, we can easily see that:
$$-(d\Gamma_{\frak{q}}+d^*Y_{\frak{q}})\wedge (d\Gamma_{\frak{q}}+d^*Y_{\frak{q}}) \in L^{2,1}(\R^4),$$
with
\begin{equation}\label{regY2bis}
\|d\Gamma_{\frak{q}}+d^*Y_{\frak{q}})\wedge (d\Gamma_{\frak{q}}+d^*Y_{\frak{q}}\|_{L^{2,1}(\R^4)}\lesssim (\| \nabla \frak{q} \|_{L^{4}}^{2} + \| \beta \|_{L^{(4,2)}})^2.
\end{equation}
From \eqref{regY2}, \eqref{nablaG} and \eqref{regY2bis} it follows 
that  $\nabla Y_{\frak{q}} \in L^{(4,1)}$ with
$$\| \nabla Y_{\frak{q}}\|_{L^{(4,1)}} \lesssim (\| \nabla \frak{q} \|_{L^{4}}^{2} + \| \beta \|_{L^{(4,2)}})^2,$$
 
\noindent
{\bf Step 3. }   We may write:
\begin{eqnarray}\label{decD7}
{\mathcal{D}}(\frak{q})&=&D\Gamma_{\frak{q}}+\psi_{\frak{q}}
\end{eqnarray}
where $\psi_{\frak{q}}\in L^{(4,1)}(\R^4)$ and $$\|\psi_{\frak{q}}\|_{L^{(4,1)}}\lesssim  (\|\nabla\frak{q}\|^2_{L^{4}} + \| \beta \|_{L^{(4,2)}})^2.$$
%{\color{red} Needs to include $\| \nabla \alpha \|_{L^{(4,2)}}$ term from the regularity of $dv$, does not affect later arguments due to $\varepsilon$-bound on this expression being chosen sufficiently small}\\
This is simply due to \eqref{hodge} and the explicit formula for $\mathcal{D}(\frak{q})$. It follows by direct evaluation of the term that
\begin{eqnarray}\label{decD2}
(\Pi_{4}+{\mathcal{P}}){\mathcal{D}}(\frak{q})=(\Pi_{4}+{\mathcal{P}})(D\Gamma_{\frak{q}})+(\Pi_{4}+{\mathcal{P}})\psi_{\frak{q}}
\end{eqnarray}
Next, we notice that $(\Pi_{4}+{\mathcal{P}}){\mathcal{D}}(\frak{q})=\beta e_4$ if and only if:
\begin{eqnarray}\label{decD3}
(\Pi_{4}+{\mathcal{P}})(D\Gamma_{\frak{q}})&=&\beta e_4-(\Pi_{4}+{\mathcal{P}})\psi_{\frak{q}}\\
&=&(\Pi_{4}+{\mathcal{P}})(\beta e_4)-(\Pi_{4}+{\mathcal{P}})\psi_{\frak{q}}.\nonumber
\end{eqnarray}
We have seen that the linear operator $(\Pi_{4}+{\mathcal{P}})\circ D$ (which in fact corresponds to the differential ${\mathcal{L}}_{\frak{I}}$ computed in the previous subsection) is an elliptic operator and if
$w=((\Pi_{4}+{\mathcal{P}})\circ D)^{-1}(\beta e_4)$ and $\tilde\psi_{\frak{q}}=-((\Pi_{4}+{\mathcal{P}})\circ D)^{-1}(\Pi_{4}+{\mathcal{P}})\psi_{\frak{q}}$, then by  Lemma \ref{4D:regP2}:

$$
\Pi_{e_{i+1}e_{i-1}e_4}(Dw)=0.
$$
From \eqref{decD2} and $w = \Gamma_{\frak{q}} - \tilde{\psi}_{\frak{q}}$,  it follows that
$$
\Pi_{e_{i+1}e_{i-1}e_4}(D\Gamma_{\frak{q}})=\Pi_{e_{i+1}e_{i-1}e_4}(D \tilde\psi_{\frak{q}})\in L^{(4,1)}(\R^4)$$
with  by elliptic estimates:
\begin{equation}
\|\nabla \tilde\psi_{\frak{q}}\|_ {L^{(4,1)}(\R^4)}\lesssim \|\nabla   Y_{\frak{q}}\|_ {L^{(4,1)}(\R^4)}.\end{equation}
  It follows now:
\begin{equation}
\Pi_{e_{i+1}e_{i-1}e_4}\left({\mathcal{D}}(\frak{q})\right)\in L^{(4,1)}(\R^4)~~~\mbox{for all $i=1,2,3.$}
\end{equation}
This shows the desired improved regularity result. \qed

\subsection{Conclusion of the Bootstrap Test}

Let $\frak{f}\in L^{4/3}(\R^4)$  be a solution of \eqref{maineq3}. By choosing $\frak{q}$ as with our gauge operator, we find:
\begin{eqnarray}\label{maineq4}
\partial_{x_0}[\frak{q}\frak{f}]-\sum_{i=1}^{3}[\partial_{x_i}(\frak{q}e_i \frak{f})]&=&
\frak{q}{V}(x)\frak{f}.
\end{eqnarray}
where ${V}(x)\in L^{(4,1)}$ by our investigation in the previous subsection (${V}(x)=\psi_q+\Pi_{e_{i+1}e_{i-1}e_4}(D\Gamma_{\frak{q}})$). Indeed, observe that this is a consequence of the choice of gauge and the improved integrability we have established. Furthermore, by the estimate proven before:
$$\|{V}(x)\|_{L^{(4,1)}}\lesssim  (\|\nabla\frak{q}\|^2_{L^{4}} + \| \beta \|_{L^{(4,2)}})^2.$$
Since 
From Lemma \ref{inj}    we can  get rid of the power 2 by choosing $\varepsilon > 0$ possibly slightly smaller. Indeed, we can show using the estimate  \eqref{Nq2} for $\frak{q}$:
$$\| V \|_{L^{(4,1)}} \lesssim \| \beta \|_{L^{(4,2)}}$$
We set 
$$F=\left(\begin{array}{c}\frak{q}\frak{f}\\-\frak{q}e_1\frak{f}\\-\frak{q}e_2\frak{f}\\-\frak{q}e_3\frak{f}\end{array}\right)$$
Our goal is to prove Morrey estimates just like in \cite{DLR1}. In order to achieve this, we will use a non-linear Hodge decomposition. The reason behind this is, that Wente's inequality is no longer at our disposal and therefore, we need a suitable replacement, see Lemma \ref{prel:compensation}.  \par
\noindent
{\bf Claim 1:} There are 
 $A,B\in \dot W^{1,(4/3,\infty)}(\R^4)$, where $B $ is differential $2$-form, such that:
\begin{equation}\label{hodgeF}
F=dA+\frak{q}d^{*}B
\end{equation}
\noindent
{\bf Proof of the Claim 1.} We argue by induction:\\

\noindent
  {\bf Step 1.} We find $A_0,B_0$ such that
  \begin{eqnarray}
-\Delta A_0&=&-\div (F)\\
-\Delta B_0&=&d(\frak{q}^{-1}F)
  \end{eqnarray}
  Then for $k\ge 1$ we solve:
    \begin{eqnarray}
   -\Delta A_k&=&-d^*(q d^*B_{k-1})=\ast (dq\wedge d\ast B_{k-1})\\
-\Delta B_k&=&-d(\frak{q}^{-1})\wedge dA_{k-1}
  \end{eqnarray}
  We set $A=\sum_{i=0}^{\infty} A_k$ and  $B=\sum_{i=0}^{\infty} B_k$.
  We then have:
  \begin{eqnarray}
   -\Delta A&=&-d^*(\frak{q} d^*B)-\div(F)\label{eqA}\\
-\Delta B&=&-d(\frak{q}^{-1})\wedge dA_{k-1}+d(\frak{q}^{-1}F).\label{eqB}
  \end{eqnarray}
From \eqref{eqA} and \eqref{eqB}, we deduce the following estimates:
  \begin{eqnarray}
 && d^{*}(F-dA-qd^{*}B)=0\label{eqdsF}\\
 &&d(\frak{q}^{-1}F-d^{*}B-\frak{q}^{-1}dA)=0\label{eqdF}
  \end{eqnarray}
  From \eqref{eqdF}, it follows there exists a function $\gamma\in \dot W^{1,4/3}(\R^4)$ such that
    \begin{eqnarray}\label{dbeta}
 &&\frak{q}^{-1}F-d^{*}B-\frak{q}^{-1}dA=d\gamma.
 \end{eqnarray}
 By combining \eqref{eqdsF} and  \eqref{dbeta} we get
 \begin{eqnarray}
 &&d^*(\frak{q}d\gamma)=0\\
 &&d(\frak{q}d\gamma)=d\frak{q}\wedge d\gamma\\
 &&\|\frak{q}d\gamma\|_{L^{(4/3, \infty)}}\lesssim\|d\frak{q}\|_{L^4}\|d\gamma\|_{L^{(4/3, \infty)}}\le \varepsilon_0 \|d\gamma\|_{L^{(4/3, \infty)}}
 \end{eqnarray}
 Notice that in the last line, we used the compensation result in Lemma \ref{prel:compensation}.\par
 It follows that, if $\varepsilon > 0$ is chosen small enough, $d\gamma=0$ and therefore
 $$
 F=dA+\frak{q}d^{*}B.$$
 We conclude the {\bf proof the claim 1}.~~\hfill$\Box$
\par
\medskip
We continue with the {\bf proof of Theorem \ref{mainth4D}}: From \eqref{hodgeF}, it follows that
\begin{equation}\label{A}
-\Delta A=\frak{q}{V}(x)\frak{f}+d^*(\frak{q} d^*B)=\frak{q}{V}(x)\frak{f}+\ast(d\frak{q} \wedge d\ast B).
\end{equation}
Then, by using the fundamental solution, we see:
\begin{eqnarray}\label{A1}
\|\nabla A\|_{L^{(4/3, \infty)}}&\lesssim&\|-\Delta A\|_{L^1}\lesssim   \|{V}\|_{L^{(4,1)}}\|\frak{q\ f}\|_{L^{(4/3, \infty)}}+\|\nabla\frak{q}\|_{L^4}\|d^*B\|_{L^{(4/3, \infty)}}\nonumber \\
&\lesssim& \|\beta \|_{L^{(4,2)}}\ \| \frak{f}\|_{L^{(4/3, \infty)}} +\|\nabla\frak{q}\|_{L^4}\|d^*B\|_{L^{(4/3, \infty)}}\nonumber\\
&\lesssim&\varepsilon \| \frak{q\,f}\|_{L^{(4/3, \infty)}}+\|\beta\|_{L^{(4,2)}}\|d^*B\|_{L^{(4/3, \infty)}}\nonumber \\
&\lesssim&\varepsilon (\| \frak{q\,f}\|_{L^{(4/3, \infty)}}+\|d^*B\|_{L^{(4/3, \infty)}})
\end{eqnarray}
Computing $\Delta B$ using exactness, we find:
 \begin{equation}\label{B}
-\Delta B=d(\frak{q}^{-1}F)+d(\frak{q}^{-1} dA)=d(\frak{q}^{-1}F)+d\frak{q}^{-1} \wedge dA
\end{equation}
From \eqref{B}, it follows as above that
\begin{eqnarray}\label{estB}
\|\nabla B\|_{L^{(4/3, \infty)}}&\lesssim&\|d\frak{q}^{-1}\|_{L^4}\|\nabla A\|_{L^{(4/3, \infty)}}+\|\frak{q}\frak{f}\|_{L^{(4/3, \infty)}}
\end{eqnarray}
By plugging \eqref{estB} into \eqref{A1}, we get for $\varepsilon > 0$ sufficiently small:
\begin{eqnarray}\label{A2}
\|\nabla A\|_{L^{(4/3, \infty)}}&\lesssim &\varepsilon \|\frak{q}\frak{f}\|_{L^{(4/3, \infty)}}
\end{eqnarray}

We set $d^*B=\sum_{i=0}^3 b_idx_i$.
By definition, it holds $d^*d^* B=\sum_{i=0}^3\partial_{x_i} b_i=0$. Moreover, by comparison of the entries in $F$, we observe: 
\begin{eqnarray*}
\frak{q} b_0&=&\frak{q} {\frak{f}}-\partial_{x_0}A\\
\frak{q} b_i&=&-\frak{q}\ e_i \frak{f}-\partial_{x_i}A 
\end{eqnarray*}
These can be slightly rearranged to express $b_j$ in terms of $\frak{f}$:
 \begin{eqnarray*}
 b_0&=& {\frak{f}}-\frak{q}^{-1}\partial_{x_0}A\\
 b_i&=&- \ e_i \frak{f}-\frak{q}^{-1}\partial_{x_i}A 
\end{eqnarray*}
 Hence, if we solve the equations above for $\frak{f}$:
  \begin{eqnarray}\label{bi}
   {\frak{f}}&=&b_0+ \frak{q}^{-1}\partial_{x_0}A=e_i (b_i+\frak{q}^{-1}\partial_{x_i}A )
  \end{eqnarray}
  Then it is now clear:
  $$\partial_{x_i}b_i=-e_i\partial_{x_i}b_0-e_i \partial_{x_{i}} \left( \frak{q}^{-1}\partial_{x_0}A \right) - \partial_{x_{i}}\left( \frak{q}^{-1}\partial_{x_i}A \right).$$
  Using the previously established fact that $\sum_{i=0}^3\partial_{x_i} b_i=0$, we note:
  \begin{eqnarray}
  \label{diracforb0}
  \partial_{x_0}b_0-\sum_{i=1}^{3} e_i\partial_{x_i}b_0=\sum_{i=1}^{3} \partial_{x_{i}} \left( e_i \frak{q}^{-1}\partial_{x_0}A+\frak{q}^{-1}\partial_{x_i}A \right) \in W^{-1, (4/3,\infty)}
  \end{eqnarray}
  As a result, using ellipticity and the corresponding estimates:
  %{\color{red} (TO CHECK that $\partial_L$ is CZ on the Clifford Algebra $\Rightarrow$ Clear as $\overline{\partial}_{L} \partial_{L} = \Delta$)}
  $$\|b_0\|_{L^{(4/3, \infty)}}\lesssim\|\nabla A\|_{L^{(4/3, \infty)}}$$
  From \eqref{bi}, this estimate easily generalises to all $b_j$. Namely, it follows that 
  $$\|b_i\|_{L^{(4/3, \infty)}}\lesssim\|\nabla A\|_{L^{(4/3, \infty)}}, ~~~\forall i=1,2,3.$$
 Consequently, recalling the definition of the $b_j$, we arrive at the desired estimate for $d^\ast B$:
     \begin{equation}
     \|d^*B\|_{L^{(4/3, \infty)}}\lesssim\|\nabla A\|_{L^{(4/3, \infty)}}
       \end{equation}
  From \eqref{hodgeF}, it finally follows that:
  $$
  \|\frak{q}\frak{f}\|_{L^{(4/3, \infty)}}\lesssim \|\nabla A\|_{L^{(4/3, \infty)}}\lesssim \varepsilon_0\|\frak{q}\frak{f}\|_{L^{(4/3, \infty)}}.
  $$
  If $\varepsilon > 0$ is chosen small enough, then $\frak{q}\frak{f}=0$ is an immediate corollary, thus establishing the bootstrap lemma.\\

\section{ The Proof of the Main Theorem \ref{regul-quatre} in 4-D}
We observe that  Theorem \ref{regul-quatre} follows similar to Theorem \ref{mainth4D}  by using localization arguments analogous to Proposition III.4 in \cite{DLR1}.
We provide here a sketch of proof in  the $4$-D case and we refer the details to \cite{Wett}.\\

First, we will briefly explain how to obtain an appropriate version of the non-linear Hodge decomposition on balls $B_r(x)$. For simplicity's sake, let us assume $x=0$, the general case is obtained by translation. Let for this $G$ be an arbitrary $1$-form in $W^{1, \frac{4}{3}}(B_r(0))$ as obtained in the proof. Then, by classical Hodge decomposition, there exist a function $A$ on $B_r(0)$ vanishing along the boundary and a $2$-form $\tilde{A}$, such that:
\begin{equation}
	d A + d^{\ast} \tilde{A} = G
\end{equation}
Next, we consider the Hodge decomposition in the same manner of $\frak{q}^{-1}d^\ast A$, again obtaining zero boundary conditions for the function $\tilde{B}$:
\begin{equation}
	d \tilde{B} + d^{\ast} B = \frak{q}^{-1} d^{\ast} \tilde{A}
\end{equation}
Thus, we have:
\begin{equation}
	G = dA + d^{\ast} A = dA + \frak{q} d^{\ast} B + \frak{q} d\tilde{B} \Rightarrow G - dA - \frak{q} d^{\ast} B = \frak{q} d\tilde{B}
\end{equation}
We observe that on $B_r(0)$:
\begin{equation}
	\Delta \tilde{B} = d^{\ast}d\tilde{B} = d^{\ast} \big{(}\frak{q}^{-1} d^{\ast} \tilde{A} \big{)} = - \ast \big{(} d \frak{q}^{-1} \wedge d (\ast \tilde{A}) \big{)}
\end{equation}
Due to the zero boundary condition, we can therefore deduce by similar arguments as in our compensation result in Lemma \ref{prel:compensation}:
\begin{equation}
	\| \nabla B \|_{L^{\frac{4}{3}}(B_r(0))} \lesssim \| d\frak{q} \|_{L^{4}(B_r(0))} \| d^{\ast} \tilde{A} \|_{L^{\frac{4}{3}}(B_r(0))} \lesssim \varepsilon \| G \|_{L^{\frac{4}{3}}(B_r(0))}
\end{equation}
So, if $\varepsilon > 0$ is sufficiently small, we can argue by iteration that there exists a solution to the non-linear Hodge decomposition as in the case of codomains of dimension 2, such that $A$ has boundary value $0$.\\

Now, to deduce local regularity, we merely have to establish slightly improved regularity and hence Morrey estimates as in \cite{DLR1}, the full regularity as in Theorem \ref{regul-quatre} follows by Morrey-bootstrapping going over to possibly smaller balls to obtain uniform powers in the Morrey estimates. Therefore, let us just point out the differences to \cite{DLR1} and our considerations in connection with the bootstrap lemma: Namely, we can estimate $A$ as in the bootstrap lemma, if we find $A,B$ for a given $B_r(x)$. More precisely, due to the boundary conditions, we will find:
\begin{align}
	\| \nabla A \|_{L^{(4/3, \infty)}(B_r(x))} 	&\lesssim \varepsilon \|\frak{q}F \|_{L^{(4/3, \infty)}(B_r(x))} + \varepsilon \| d^{\ast} B \|_{L^{(4/3, \infty)}(B_r(x))} \notag \\
										&\lesssim \varepsilon \|\frak{q}F \|_{L^{(4/3, \infty)}(B_r(x))} + \varepsilon \| \nabla A \|_{L^{(4/3, \infty)}(B_r(x))},
\end{align}
by using the same arguments as before and using $F = dA + \frak{q} d^{\ast} B$. So if $\varepsilon$ is sufficiently small, we arrive at:
\begin{equation}
	\| \nabla A \|_{L^{(4/3, \infty)}(B_r(x))} \lesssim \varepsilon \| \frak{q}F \|_{L^{(4/3, \infty)}(B_r(x))},
\end{equation}
Then, it remains to obtain appropriate estimates for $d^{\ast}B$. For this, write $d^\ast B = \sum_{j} b_j dx_j$ and we can deduce completely analogous to \eqref{diracforb0} in the proof of the bootstrap lemma:
$$\partial_L b_0 = \sum_{j \geq 2} \partial_{x_j} R_j,$$
where $R_j$ is an expression depending on $\frak{q}$ and $\nabla A$ as already found in the proof of Theorem \ref{mainth4D}. So we can now split $b_0$ into a Clifford analytic and thus harmonic part, which can be estimated by means of Campanato-estimates as in \cite{DLR1} and the convolution of the RHS in the equation above with the fundamental solution of $\partial_L$ on $\R^4$. This second summand can be estimated by usual estimates for the fundamental solution of the Laplacian. Therefore, we arrive at the desired estimates for $d^{\ast}B$ by completely the same means as in \cite{DLR1} once we use the link between $b_j$ and $b_0$ established in the bootstrap lemma. For details, see \cite{Wett}.

\section{ The 3-D Case}
  Before we briefly discuss the general case, let us provide another example on how to construct an appropriate gauge operator. More precisely, we shall consider the case of $3$D-domains. This will illustrate that the result we have obtained will not generalise in an "easy" manner to arbitrary dimensions $m \geq 3$, but one has to take some care when investigating the gauge operators involved:\\

\noindent
  Let us consider the following equation:
  \begin{equation}
  \label{3dmain}
  	\partial_{L} \frak{f} = \beta e_3 \cdot \frak{f},
  \end{equation}
  where $\frak{f}: \R^3 \to {C\ell}_{3}$ is in $L^{3/2}$. Let us assume that 
  $$\beta = \beta^{0} + \beta^{1} e_1 + \beta^2 e_2 \in L^{(3,2)}(\R^3;\operatorname{span}_{\R} \{ e_0, e_1, e_2 \}),$$
  as well as:
  \begin{equation}
  	\mathbf{curl}_{x_1,x_2} \beta = \partial_{2} \beta^{1} - \partial_{1} \beta^{2} = 0.
  \end{equation}
  We will sketch the proof of the following Theorem which is along the same lines as the proof of Theorem \ref{regul-quatre}:
   \begin{Th}\label{regul-trois}
 Let    $\beta=(\beta_0,\beta_1,\beta_2)\in W^{1,2}({\R}^3,\operatorname{span}_{\R}\{e_0,e_1,e_2\})$ with 
 \be
 \label{curl-03d}
  \p_{x_2}\beta_1-\p_{x_1}\beta_2=0\ .
 \ee   
 Let $\frak{f}\in L^{3/2}({\R}^3,{C\ell}_{2}^2)$ be a solution of
 \be
 \label{main-equation3d}
 \partial_L\frak{f}= \left(
 \begin{array}{cc}
 0 &\beta\\[3mm]
 -\beta & 0
 \end{array}
 \right)\ \hat{\frak{f}}
 \ee
  Then  $\frak{f}\in L^ {q}_{loc}(\R^3)$ for all $q<\infty$.   \hfill $\Box$
 \end{Th}
 It is clear that we may reformulate \eqref{main-equation3d} into an equation of the following form:
 $$\partial_L \frak{g}=\beta  e_3 \cdot \frak{g},$$
 for $\frak{g} = \frak{f}^1 + \frak{f}^2 e_3$.
 Moreover, there is also the following bootstrap test:
 \begin{Th}\label{mainth3D} There exists $\varepsilon_0>0$ such that for every $\beta\in L^{(3,2)}(\R^3, \operatorname{span}_{\R} \{ e_0, e_1, e_2 \} )$ satisfying $\| \beta\|_{L^{(3,2)}(\R^3)}\le \varepsilon_0$ as well as: 
 $$\p_{2} \beta^1 - \p_1 \beta^2 = 0,$$
 and every  $\frak{f} \in L^{3/2}(\R^3,{C\ell}_3)$ solving:
 \begin{equation}\label{maineq63d}
\partial_L \frak{f}=\beta  e_3 \cdot \frak{f}\ ,
\end{equation}
we have ${\frak{f}}\equiv 0.$ \hfill $\Box$
\end{Th}

In our current discussion, we focus on Theorem \ref{mainth3D}, see the discussion in the previous section regarding the proof of Theorem \ref{mainth4D} for a sketch on how to apply Morrey-estimates and \cite{Wett}. For later convenience, let us introduce the following spaces:
  \begin{align}
  	V_3 		&:= \operatorname{span}_{\R} \{ e_3, e_1 e_3, e_2 e_3 \} \\
	V_2 		&:= \operatorname{span}_{\R} \{ e_1, e_2 \} \\
	V_1	 	&:= \R \cdot e_1 e_2,
  \end{align}
  and denote by $\Pi_3, \Pi_2$ and $\Pi_1$ the projections of $\mathcal{U}_{3}$ onto the respective subspaces.
  
  As in \cite{DLR1} and previously seen in the case of $4$-dimensional domains, let us multiply both sides of \eqref{3dmain} by a function $\frak{q}: \R^3 \to Spin(4)$ to reveal a slight gain in integrability after a change of gauge. We obtain by using Leibniz' rule:
 \begin{eqnarray}\label{3dqf}
 \frak{q}\partial_L \frak{f}=\partial_{x_0}(\frak{q} f)-(\partial_{x_0}\frak{q})\frak{f}-\sum_{i=1}^{2}\partial_{x_i}(\frak{q} e_i\frak{f})+\sum_{i=1}^{2}\partial_{x_i}\frak{q} e_i \frak{f}
 \end{eqnarray}
 We denote by:
 $$
 {\mathcal{D}}(\frak{q}):=\frak{q}^{-1}\partial_{x_0}\frak{q} - \sum_{i=1}^{2}\frak{q}^{-1}\partial_{x_i}\frak{q} e_i = \frak{q}^{-1} \partial_R \frak{q}.$$
   Observe that
 \begin{equation}
 \label{decalpha3d}
 \beta e_3=\beta^0 \cdot e_3-\sum_{i=1}^2 \beta^{i} \cdot e_ie_3 \in V_3\\
 \end{equation}
By using \eqref{3dqf} and rearranging, we get:
 \begin{eqnarray}\label{3dneweq}
 \partial_{x_0}(\frak{q} f)-\sum_{i=1}^{2}\partial_{x_i}(\frak{q} e_i\frak{f})&=&\frak{q}(\beta\,e_3+{\mathcal{D}}(\frak{q})) \frak{f}.
 \end{eqnarray}
 We notice that in \eqref{3dneweq}, the absorption of $\beta e_3$ by $\mathcal{D}(\frak{q})$ leads to a system of  $8$ equations in merely $6$ unknowns, which is overdetermined much like in the $4$-dimensional case. Therefore, there is generally no hope of completely absorbing the "bad term", however, inspired by our proof in $4$D, we hope to absorb $\beta e_3$ up to a term of higher integrability as before.\\
 
 \noindent
   The main aim is to find $\frak{q}\in \dot W^{1,3}(\R^3, Spin(4))$ such that
   ${\mathcal{D}}(\frak{q})=-\beta e_3+V(x)$ where $V$ is a more regular potential than $\beta e_3$, namely  $V\in L^{(3,1)}(\R^3)$. To do this, let us introduce the following non-linear operator reminiscent of \eqref{gaugeop}:
  \begin{eqnarray}\label{3dest9}
&& {\mathcal{N}} \colon \dot W^{1,3}(\R^3 ,Spin(4))\to  W^{-1,3}(\R^3 ,V_2)\times   L^3(\R^3,\mathcal{U}_2) \\[5mm]
&&~~~~~ \frak{q}\mapsto\left(\Pi_2\left(\sum_{i=0}^3(\partial_{x_i}(\frak{q}^{-1}\partial_{x_i} \frak{q})\right), -\Pi_{3}({\mathcal{D}}(\frak{q})) e_3 + \Pi_1(\frak{q}^{-1} \partial_{x_{0}} \frak{q}) - \sum_{j=1}^2\Pi_1(\frak{q}^{-1} \partial_{x_{j}} \frak{q}) e_j \right) \nonumber
\end{eqnarray}
We notice that the first component is analogous to \eqref{gaugeop}, while the second component of $\mathcal{N}$ looks more complicated than before. As we shall see later, this definition neatly connects the differential of $\mathcal{N}$ to the Riemann-Fueter operator once again. Indeed, analogous to our previous discussion for $\R^4$, we have the following as a main result:
  
  \begin{Lm}\label{3dinj}
   There exists $\varepsilon_0>0$ and $C>0$ such that for any choice $\omega \in W^{-1,3}(\R^3, V_2)$ and $\frak{g} \in L^{3}(\R^3, {C\ell}_{2})$ satisfying
   \begin{equation}
   	\| \omega \|_{W^{-1,3}} \leq \varepsilon_{0}, \| \frak{g} \|_{L^{3}} \leq \varepsilon_{0},
   \end{equation}
   there is $\frak{q}\in \dot W^{1,3}(\R^3 ,Spin(4))$ such that 
   \begin{eqnarray}
   \label{3ddecomp}
   	\mathcal{N}(\frak{q}) = (\omega, \frak{g})
   \end{eqnarray}
   as well as
   \begin{equation}
   \|\nabla\frak{q}\|_{L^3} \leq C ( \| \omega \|_{W^{-1,3}} + \| \frak{g} \|_{L^{3}} ).
   \end{equation}
 \end{Lm}
The proof essentially proceeds as in \cite{DLR1} and the case of domains of dimension $4$, so let us introduce the analogous simplifications: Again similar to \cite{DLR2,DLS2}, using an approximation argument similar to the our closedness argument later on, it suffices to prove Lemma~\ref{3dinj} for $\om$ and $\frak{g}$ slightly more integrable, namely under the assumption $\om\in ( W^{-1,p}\cap  W^{-1,p'})({\R^3},V_2)$ and 
 $\frak{g}\in   (L^p\cap L^{p^{\prime}})({\R^3},\mathcal{U}_2)$  for some $3<p$, $p'=\frac{p}{p-1}$. For the remainder of our discussion, we fix some $3<p$. Given $\varepsilon > 0$, we again define as previously:
\be
\label{3dUepsilon}
     C\ell_{\varepsilon}  := \left \{
     \begin{array}{c}
     (\omega,\frak{g}) \in  (W^{-1,p}\cap   W^{-1,p'})({\R^2},V_2)\times (L^p\cap L^{p^{\prime}})({\R^3},\mathcal{U}_2)\\[5mm]
   \|\omega\|_{W^{-1,3}} +\|\frak{g}\|_{L^3}\le \varepsilon
  \end{array}
  \right \}
  \ee

For constants $\varepsilon, \Theta > 0$, let $\mathcal{V}_{\varepsilon,\Theta} \subseteq C\ell_{\varepsilon}$ be the set where we have the decomposition \eqref{3ddecomp} with the estimates
 \be
\label{3deq:gauge:2est}
  \|\nabla  \frak{q}\|_{L^3} \leq \Theta  ( \|\omega\|_{W^{-1,3}}+\|\frak{g}\|_{L^3})\, 
 \ee 
 \be
  \label{3deq:gauge:pest}
 \|\nabla \frak{q}\|_{p} \leq \Theta  ( \|\omega\|_{W^{-1,p}}+\| \frak{g}\|_{L^p})\, ,
 \ee
 \be
 \label{3deq:gauge:ppest}
     \|\nabla \frak{q}\|_{p^{\prime}} \leq \Theta  ( \|\omega\|_{W^{-1,p^\prime}}+\|\frak{g}\|_{L^{p^\prime}})\,.
\ee
That is
\[
 \mathcal{V}_{\varepsilon,\Theta} := \left \{\omega,\frak{g} \in C\ell_{\varepsilon}:\ \begin{array}{c}
                                             \mbox{there exists $\frak{q} \in (\dot{W}^{1,p} \cap \dot W^{1,p^{\prime}})(\R^3,Spin(4))$, so that  }\\[3mm]
                                           \frak{q}-\frak{I}\in  L^{3p/2p-3} (\R^3,Spin(4))\\[3mm]~~\mbox{and}~~ \eqref{3ddecomp}, \eqref{3deq:gauge:2est}, \eqref{3deq:gauge:pest},  \eqref{3deq:gauge:ppest}~~ \mbox{hold.}
                                              \end{array}\right \}
\]
The strategy   to prove Lemma~\ref{3dinj} is precisely the same as for Lemma \ref{lm-invertN-a} and it is a corollary of the following:
\begin{Prop}\label{3dpr:VepseqUeps}
There exist $\Theta > 0$ and $\eps > 0$, such that $\mathcal{V}_{\eps,\Theta} = C\ell_{\varepsilon}$. \hfill $\Box$
\end{Prop}
\noindent{\bf Proof of Proposition \ref{3dpr:VepseqUeps}.}
Proposition~\ref{3dpr:VepseqUeps} follows, once we show the following four properties
\begin{itemize}
 \item[(i.)] $C\ell_{\varepsilon}$ is  connected.
 \item[(ii.)] $\mathcal{V}_{\eps,\Theta}$ is nonempty. 
 \item[(iii.)] For any $\eps, \Theta > 0$, $\mathcal{V}_{\eps,\Theta}$ is a relatively closed subset of $C\ell_{\varepsilon}$.
 \item[(iv.)] There exist $\Theta > 0$ and $\eps > 0$ so that $\mathcal{V}_{\eps,\Theta}$ is a relatively open subset of $C\ell_{\varepsilon}$.
\end{itemize}
\par
As in \cite{DLR1}, property (i.) and (ii.) are obvious and (iii.) follows as in the case of 4-dimensional domains. For further details, we refer to our discussion of the $4$D-case.
\medskip  

It remains to show the {openness property} (iv.). For this let $(\omega_0,\frak{g}_0)$ be arbitrary in $\mathcal{V}_{\eps,\Theta}$. Let $\frak{q}_0 \in \dot{W}^{1,p} \cap \dot{W}^{1,p^{\prime}}(\R^3,Spin(4))$, $ \frak{q}_0-\frak{I} \in L^{3p/2p-3}(\R^3 )$ so that the decomposition \eqref{3ddecomp} as well as the estimates \eqref{3deq:gauge:2est}, \eqref{3deq:gauge:pest} and \eqref{3deq:gauge:ppest}  are satisfied for $\omega_0$ and $\frak{g}_0 $. The idea is to study perturbations of $\frak{q}_0$ of the form $\frak{q}=\frak{q}_{0}e^{\frak{u}}$, where $\frak{u}\in \dot W^{1,p} \cap \dot W^{1,p^{\prime}}(\R^3,\frak{spin}(4)) \cap L^{3p/2p-3}(\R^{3})$. Completely analogous to before, the exponent $p>3$ has been chosen in particular to ensure $p^\prime < 3$ and, as a result, $\frak{u}\in C^0\cap L^\infty({\R^3})$ and $\frak{q}_{0}e^{\frak{u}}-\frak{I}\in L^{\frac{3p}{2p-3}}$. This follows precisely the same way as in the 4-dimensional case treated previously, where we mentioned that the main estimate is independent of the dimension of the underlying space.\par

Attentive readers know what comes next: We compute the differential $D {\mathcal{N}}  (\frak{q}_0)$ as 
\[
D {\mathcal{N}}  (\frak{q}_0)= \frac{d}{dt} {\mathcal{N}} (\frak{q}_0e^{t\frak{u}}) \Big|_{t=0}=: {\mathcal{L}}_{\frak{q}_{0}}({\frak{u}}),
\]
where $\frak{u}\in (\dot W^{1,p} \cap \dot W^{1,p^{\prime}} \cap L^{3p/2p-3})(\R^3,\frak{spin}(4))$. We write
$${\mathcal{L}}_{\frak{q}_{0}}({\frak{u}})=({\mathcal{L}}^2_{\frak{q}_{0}}({\frak{u}}),{\mathcal{L}}^{3}_{\frak{q}_{0}}({\frak{u}}))$$
where 
\begin{eqnarray*}
{\mathcal{L}}^2_{\frak{q}_{0}}({\frak{u}})&:=&\Pi_{2}\left[\Delta\frak{u}+
\sum_{j=0}^2\partial_{x_j}\left(\frak{q}_0^{-1}(\partial_{x_j}\frak{q}_0)\frak u-\frak{u}\frak{q}^{-1}_0\partial_{x_j}\frak{q}_0\right)\right]\\
{\mathcal{L}}^{3}_{\frak{q}_{0}}({\frak{u}})&=&-\Pi_{3}(\partial_{R} \frak{u})e_3 + \partial_{R} \Pi_{1}(\frak{u})\nonumber\\[3mm]&&-\sum_{j=0}^2(-1)^{\delta_{0j}}\Pi_{3}\left( (\frak{q}_{0}^{-1} \partial_{x_{j}}\frak{q}_{0} \frak{u} - \frak{u} \frak{q}_{0}^{-1} \partial_{x_{j}}\frak{q}_{0})e_{j} \right) e_3\nonumber\\[3mm]&&+ \sum_{j=0}^{2} (-1)^{\delta_{0j}} \Pi_{1}(\frak{q}_{0}^{-1} \partial_{x_{j}}\frak{q}_{0} \frak{u} - \frak{u} \frak{q}_{0}^{-1} \partial_{x_{j}}\frak{q}_{0})e_j
\end{eqnarray*}
The essential property we will be using is the invertibility of ${\mathcal{L}}_{\frak{q}_0}({\frak{u}})$ in the special case $\frak{q}_0=\frak{I}$.
If $\frak{q}_0=\frak{I}$, we have $d\frak{q}_{0} = 0$ and therefore the differential simplifies significantly:
\begin{eqnarray}
{\mathcal{L}}^2_{\frak{I}}({\frak{u}})&=&\Pi_{2}\left[\Delta\frak{u}\right]\nonumber\\
{\mathcal{L}}^{3}_{\frak{I}}({\frak{u}})&=&-\Pi_{3}(\partial_{R} \frak{u})e_3 + \partial_{R} \Pi_{1}(\frak{u})
 \end{eqnarray}

\medskip 
\begin{Prop}\label{3dinvert}
 The operator ${\mathcal{L}}_{\frak{I}}({\frak{u}})$ is elliptic. \end{Prop}
\par
\noindent
{\bf Proof of Proposition \ref{3dinvert}.}
We write $\frak{u}=w+v$ where $w\in  V_2$ and $v=v^0e_3+v^1e_1e_3+v^2e_2e_3+v^3e_1e_2\in V_1 \oplus V_3$.
We observe that 
\begin{eqnarray}
{\mathcal{L}}^2_{\frak{I}}({\frak{u}})&=&\Pi_{2}\left[\Delta w + \Delta v\right] = \Delta w \nonumber\\
{\mathcal{L}}^{3}_{\frak{I}}({\frak{u}})&=&-\Pi_{3}(\partial_{R} v)e_3 + \partial_{R} \Pi_{1}(v)\nonumber
\end{eqnarray}
Computing ${\mathcal{L}}^{3}_{\frak{I}}({\frak{u}})$ explicitly, we find:
 \begin{eqnarray}
-\Pi_{3}(\partial_{R} \frak{u})e_3 + \partial_{R} \Pi_{1}(\frak{u})&=&
(\partial_{x_{0}} v^0 - \partial_{x_1}v^1 - \partial_{x_2} v^2) + (\partial_{x_1} v^0 + \partial_{x_0} v^1 - \partial_{x_2} v^3) e_1\nonumber\\[3mm]&+& (\partial_{x_2} v^0 + \partial_{x_0} v^2 + \partial_{x_1} v^3) e_2 + (\partial_{x_2} v^1 - \partial_{x_1} v^2 + \partial_{x_0} v^3) e_1e_2\nonumber\\[3mm]&=& D^{RF}_{R}(v^0 + v^1 i + v^2 j + v^3 k)\nonumber
\end{eqnarray}
We can associate to this operator the following symbol:
\begin{equation}
\sigma(\xi)=\left(\begin{array}{cccc}
\xi_0 & -\xi_1&-\xi_2&0\\[3mm]
\xi_1&\xi_0&0&-\xi_2\\[3mm]
\xi_2&0&\xi_0&\xi_1\\[3mm]
0&\xi_2&-\xi_1&\xi_0\end{array}\right)\end{equation}
It is immediately clear that this is now the Riemann-Fueter operator applied to functions depending only on the first $3$ variables. Therefore, one may argue as in $4$D that the symbol is everywhere invertible. In fact, this is an immediate corollary of the computations in $4$D. This concludes the proof of Proposition \ref{3dinvert}.~~\hfill$\Box$\par
\medskip
 We can prove the following result, which we only state, since the proof is now more or less a copy of the corresponding result in $4$D:
\begin{Lm}\label{la:invertible}
For any $\Theta > 0$, there exists  $\eps > 0$ so that the following holds for any $\omega_0,\frak{g}_0 $  and $\frak{q}_0$ as above:
\par
For any $\om\in ( W^{-1,p}\cap  W^{-1,p'})({\R^3},V_2)$ and 
 $\frak{g}\in   (L^p\cap L^{p^{\prime}})({\R^3},\mathcal{U}_2)$ 
 there exists a unique ${\frak{u}}\in   \dot W^{1,p} \cap \dot W^{1,p^{\prime}}\cap L^{3p/2p - 3}(\R^3,\frak{spin}(4))$  so that 
\[
(\omega,\frak{g}) = \mathcal{L}_{\frak{q}_0}({\frak{u}}) 
\]
and for some constant $C = C(\omega_0,\frak{g}_0,\Theta) > 0$ it holds
\begin{eqnarray}\label{esteta}
 \|\nabla{\frak{u}}\|_{L^{p}(\R^3)} +\|\nabla{\frak{u}}\|_{L^{p^{\prime}}(\R^3)} &\lesssim&
  \|\omega\|_{W^{-1,p}(\R^3)} +\|\omega\|_{W^{-1,p'}(\R^3)}\\&&+\|\frak{g}\|_{L^{p}(\R^3)} +\| \frak{g}\|_{L^{p^{\prime}}(\R^3)}.\nonumber 
\end{eqnarray}
\end{Lm}

\noindent{\bf Proof of Proposition \ref{3dpr:VepseqUeps} continued.}

  \par
  \noindent
 For $\eps = \eps(\Theta) > 0$ chosen small enough and for any $(\omega_0,\frak{g}_0) \in \mathcal{V}_{\eps,\Theta}$, the local inversion theorem applied to ${\mathcal{N}}_{\frak{q}_0}$ gives the existence of some $\delta > 0$ (that might depend on $(\omega_0,\frak{g}_0)$) such that, for every $(\omega,\frak{g}) \in C\ell_{\varepsilon}$ with

\begin{eqnarray}
&& \|\omega - \omega_0 \|_{W^{-1,p}(\R^3)}+\|\omega - \omega_0 \|_{W^{-1,p^{\prime}}(\R^3)} < \delta\label{3destgauge5}\\[5mm]
&& \|\frak{g} - \frak{g}_0 \|_{L^{p}(\R^3)}+\|\frak{g}- \frak{g}_0 \|_{L^{p^{\prime}}(\R^3)} < \delta, \label{3destgauge6}
\end{eqnarray}
we find $\frak{q} = \frak{q}_0 e^{\frak{u}} \in \dot{W}^{1,p} \cap \dot{W}^{1,p^{\prime}}(\R^3, Spin(4))$, so that $ \frak{q}-\frak{I}  \in  L^{3p/2p-3}(\R^3)$ and \eqref{3ddecomp} is satisfied.
It remains to prove \eqref{3deq:gauge:2est}, \eqref{3deq:gauge:pest} and \eqref{3deq:gauge:ppest}. This will be implied by the following lemma, whose proof is again analogous to the $4$D-case and therefore omitted:
\begin{Lm}\label{3dla:qgaugeuniformest}
There exists $\Theta > 0$ and $\sigma > 0$, such that whenever $\frak{q} \in \dot{W}^{1,p} \cap \dot{W}^{1,p^{\prime}}(\R^3)$  with $ \frak{q}-\frak{I}  \in  L^{3p/2p-3}(\R^3)$ satisfying \eqref{3ddecomp} and it holds
\begin{equation}
  \|\nabla\frak{q}\|_{L^3(\R^3)}\leq \sigma,
\end{equation}
then  \eqref{3deq:gauge:2est}, \eqref{3deq:gauge:pest}
                                            and \eqref{3deq:gauge:ppest} 
 hold true as well.\hfill $\Box$
\end{Lm}

Thanks to Lemma~\ref{3dla:qgaugeuniformest}, the openness property (iv.) is proven. Proposition~\ref{3dpr:VepseqUeps}, and as a corollary also Lemma \ref{3dinj}, is now established.
\hfill $\Box$\\

In order to establish the bootstrap lemma and Morrey estimates, one can now proceed completely analogous to the case of domains of dimension $4$. Indeed, the arguments for improved regularity of the potential carry over immediately and the non-linear Hodge decomposition works equally well in this case. We refer to our discussion for $\R^4$, the modifications should be self-explanatory.

\section{ Perspectives for Domains of Dimension $\mathbf{ 5 \leq m \leq 8}$}

Finally, let us briefly discuss the possibility to extend the results presented to domains of arbitrary dimensions $\leq 8$. There in fact is a way to generalise the construction of the gauge operator in these cases and we refer to \cite{Wett} for the details. The key is that in the cases $m=3$ and $m=4$, the gauge operator relies on the ellipticity of the Riemann-Fueter operator to show existence and approrpriate estimates. For $5 \leq m \leq 8$, we may substitute the Riemann-Fueter operator by the octonionic derivative in a suitable sense, which allows us to conclude in much the same way. This is not very surprising, considering that the Riemann-Fueter operator is indeed the same as the quaternionic derivative. In some sense, the main property we use is the existence of an orthogonal frame which happens to parallelize the sphere, a property closely linked to the existence of normed division algebras and thus to quaternions and octonions. Since this is only possible for the spheres in dimension $0,1,3,7$, we are thus restricted by our technique to $m \leq 8$. If one manages to find a sufficiently nice elliptic, first order operator having some additional properties to ensure that it is related to the change of gauge as in \eqref{estaux}, the range of dimensions $m$ to which our proof applies could be extended.

\appendix

\section{Riemann-Fueter and Dirac operators}\label{RF}
In this appendix, we introduce and define the most important notions  that have been used in this note. We mostly limit ourselves to stating the definitions and main properties and refer to the literature for further details as well as the corresponding proofs.

The reduction from a system of divergence PDE to a linear one will be greatly simplified by introducing a family of important first order differential operators, the so-called Dirac operators. In one of the final sections, we shall consider a variation of the definition here which retains most of the same properties, but is slightly better behaved with respect to the change of gauge we envision.

\subsection{Riemann-Fueter Operator on $\K$}
Let $f: \K \to \K$ be a quaternion-valued function over $\K \simeq \R^4$. The $4$D-Riemann-Fueter operator $D^{RF}_L$ acting from the right is defined by:
\begin{align}
	D^{RF}_{R} f := 	&\big{(} \partial_{x_{0}} f_{0} - \partial_{x_{1}} f_{1} - \partial_{x_{2}} f_{2} - \partial_{x_{3}} f_{3} \big{)} \notag \\
				&+ \big{(} \partial_{x_{0}} f_{1} + \partial_{x_{1}} f_{0} - \partial_{x_{2}} f_{3} + \partial_{x_{3}} f_{2} \big{)}i \notag \\
				&+ \big{(} \partial_{x_{0}} f_{2} + \partial_{x_{1}} f_{3} + \partial_{x_{2}} f_{0} - \partial_{x_{3}} f_{1} \big{)}j \notag \\
				&+ \big{(} \partial_{x_{0}} f_{3} - \partial_{x_{1}} f_{2} + \partial_{x_{2}} f_{1} + \partial_{x_{3}} f_{0} \big{)} k,
\end{align}
where $f = f_{0} + f_{1} \cdot i + f_{2} \cdot j + f_{3} \cdot k$, or abbreviated:
$$D^{RF}_{R} f = \partial_{x_0} f + \partial_{x_1} f \cdot i + \partial_{x_2} f \cdot j+ \partial_{x_3} f \cdot k.$$
The conjugated differential operator $\overline{D}^{RF}_R$ is similarily defined:
$$\overline{D}^{RF}_{R} f = \partial_{x_0} f - \partial_{x_1} f \cdot i - \partial_{x_2} f \cdot j - \partial_{x_3} f \cdot k.$$
It is easy to see by a direct calculation:
$$\overline{D}^{RF}_{R} D^{RF}_{R} f = D^{RF}_{R} \overline{D}^{RF}_{R} f = \Delta f.$$
This can for instance be proven by considering the symbol $\sigma_{D^{RF}_{R}}$ of the differential operator $D^{RF}_{R}$:
\begin{equation}
\sigma_{D^{RF}_{R}}(\xi)=\left(\begin{array}{cccc}
\xi_0 & -\xi_1&-\xi_2&-\xi_3\\[3mm]
\xi_1&\xi_0&\xi_3&-\xi_2\\[3mm]
\xi_2&-\xi_3&\xi_0&\xi_1\\[3mm]
\xi_3&\xi_2&-\xi_1&\xi_0\end{array}\right)\end{equation}
We emphasize that the connection between $D^{RF}_{R}$ and the Laplacian mirrors the same relation between the complex derivative $\partial_{z}$ and the Laplacian. In particular, we have access to regularity results by using the Laplacian as an intermediate step. In particular, deriving a fundamental solution is greatly simplified and many results from complex analysis can be carried over to Riemann-Fueter operators, see \cite{garling}. As a simple example, if $D^{RF}_{R} f = 0$, then $f$ is automatically harmonic and thus smooth.\\

Naturally, analogous operators $D^{RF}_{L}$ and $\overline{D}^{RF}_{L}$ using multiplication from the left rather than from the right can be defined and satisfies similar properties. However, it should be noted, that the two pairs of operators are not the same due to the non-commutativity of the quaternions. This is in stark contrast with the situation on $\C$, which is a commutative field, and already hints at possible difficulties that might arise in our arguments later on.

\subsection{General Dirac Operators on Clifford Algebras}
Let now $m \in \N$ be given and we define for functions $f: U \subset \R^{m+1} \to {C\ell}_{m}$ the Dirac operator $\partial_{L}$ in the following way:
\begin{equation}
	\partial_{L} f = \partial_{x_{0}} f - e_{1} \cdot \partial_{x_{1}} f - \ldots - e_{m} \cdot \partial_{m} f.
\end{equation}
We refer to \cite{garling} for details on properties of this kind of operator. Once again, we can easily generalise this definition by changing signs to obtain $\overline{\partial}_{L}$ or by moving the multiplications to the other side to arrive at $\partial_{R}$ and $\overline{\partial}_{R}$ respectively.\\

By a direct computation, we can easily deduce that:
$$\partial_{L} \overline{\partial}_{L} f = \overline{\partial}_{L} \partial_{L} f = \Delta f,$$
extending the connection between the Laplacian and complex differentiation or the Riemann-Fueter operator to arbitrary Clifford algebras. We emphasise that the Riemann-Fueter operator is not a special case of the Dirac operators, although they share a lot of common features, see \cite{garling}. In addition, observe the different conventions regarding the signs associated with the partial derivatives. As earlier, this enables us to easily extend regularity results for the Laplacian to the Dirac operators.\\

For completeness' sake, let us introduce the following notion as in \cite{garling}: A function $f$ is called \emph{Clifford-analytic}, if $\partial_{L} f = 0$. By our previous elaborations, such functions are harmonic and thus smooth. A theory of such functions in analogy to complex analysis can be built up from scratch, see \cite{garling} as well as the theory of Hardy spaces by using Clifford analytic functions.

\subsection{Spin Groups}
An important subset of ${C\ell}_m$ is the so-called \emph{Spin-group}: For a fixed $m \in \N$, we define:
$$\mbox{Spin}(m) = \{ v_{1} \cdot \ldots \cdot v_{2k}\ |\ k \in \N, v_{j} \in {C\ell}^{(1)}_m \simeq \R^m \text{ and } \| v_{j} \| = 1 \text{ for all } j \} \subset {C\ell}_{m}$$
These groups are actually compact Lie groups and provide a natural two-fold covering of $\mathfrak{so}(m)$. Their Lie algebras are given by:
$$\frak{spin}(m) = {C\ell}^{(2)}_m.$$
Observe that $\operatorname{dim}_{\R} \mbox{Spin}(m) = \frac{1}{2} m(m-1)$. In a similar manner, we can introduce the compact Lie groups $\mbox{Spoin}(m)$, see \cite{garling}:
$$\mbox{Spoin}(m) = \{ v_1 \cdot \ldots \cdot v_{k}\ |\ v_j \in {C\ell}^{(0)}_m \oplus {C\ell}^{(1)}_m \simeq \R^{m+1} \text{ and } \| v_{j} \| = 1 \text{ for all } j \} \subset {C\ell}_{m}$$
This group provides another two-fold covering, this time one for $\mathfrak{so}(m+1)$. As a result, it is easy to deduce that $\mbox{Spoin}(m) \simeq \mbox{Spin}(m+1)$ due to the uniqueness of the universal covering of $\mathfrak{so}(m+1)$. The Lie algebra $spoin(m)$ is given by:
$$spoin(m) = {C\ell}^{(1)}_m \oplus {C\ell}^{(2)}_m \simeq \frak{spin}(m+1).$$
As a simple, explicit example, we have:
$$\frak{spin}(4) \simeq spoin(3) \simeq \operatorname{span} \{ e_{1}, e_2, e_{3}, e_1 e_2, e_1 e_3, e_2 e_3 \}.$$
In what follows, we will usually denote $Spoin(m)$ in ${C\ell}_{m}$ by $Spin(m+1)$ in order to adhere to common terminology. We refer to Theorems 6.3, 6.8, 6.12, 7.26, 7.27 and 8.10 in \cite{garling} for further details regarding these groups.

\subsection{Hodge Decomposition and Hodge $\ast$-Operator}
Let us briefly recall the Hodge $\ast$-operator on $\R^m$ with respect to the standard basis. On $\R^m$, we use the standard basis $b_0, \ldots, b_{m-1}$ and we have for the standard euclidean inner product:
$$\langle b_i, b_j \rangle = \delta_{ij}, \quad \forall i,j \in \{0, \ldots, m-1\}.$$
Denote by $b_0^*, \ldots, b_{m-1}^*$ the dual basis. Then $b_{i_1}^* \wedge \ldots \wedge b_{i_{k}}^*$ for $0 \leq k \leq m$ and $0\leq i_1 < \ldots < i_k \leq m-1$ form a basis for $\bigwedge \R^m$. We may now define a scalar product $\langle \cdot, \cdot \rangle_{\bigwedge \R^m}$ on $\bigwedge \R^m$ by declaring the collection of all $b_{i_1}^* \wedge \ldots \wedge b_{i_{k}}^*$ to be an orthonormal basis. The scalar product can also be defined, and actually is, independent of the choice of orthonormal basis $b_0, \ldots, b_{m-1}$, even for arbitrary $k$-forms as well as arbitrary Riemannian metrics $g$, by using local $g$-orthonormal frames. From now on, we shall write $dx_j$ instead of $b_j^*$, following the usual convention.\par
The Hodge $\ast$-operator is then defined for all $\eta, \omega$ $k$-forms by the following formula:
$$\eta \wedge \ast \omega = \langle \eta, \omega \rangle_{\bigwedge \R^m} \mu,$$
where $\mu = dx_0 \wedge \ldots \wedge dx_{m-1}$ is the standard volume form on $\R^m$. Using this operator, we can introduce the codifferential $d^\ast$ of a $k$-differential form $\omega$ on $\R^m$ by the following formula:
$$d^\ast \omega = (-1)^{m(k-1) + 1} \ast d \ast \omega,$$
where $d$ denotes the usual exterior derivative on differential forms. The Laplacian of a form $\omega$ is then defined as follows:
$$-\Delta \omega = (dd^\ast + d^\ast d) \omega.$$
Let us provide a computation of $d^{\ast}$ in the special case $m=4$: Assume $\omega = \omega^0 dx_0 + \ldots + \omega^3 dx_3$ is a $1$-form. Direct considerations show that:
\begin{align*}
	d^\ast \omega 	&= - \star d \big{(} \omega^0 dx_1 \wedge dx_2 \wedge dx_3 -+ \ldots - \omega^3 dx_0 \wedge dx_1 \wedge dx_2 \big{)}\\
				& = - \star \big{(} \partial_{x_0} \omega^0 + \partial_{x_1} \omega^1 + \partial_{x_2} \omega^2 + \partial_{x_3} \omega^3 \big{)} \mu = - \big{(} \partial_{x_0} \omega^0 + \partial_{x_1} \omega^1 + \partial_{x_2} \omega^2 + \partial_{x_3} \omega^3 \big{)}
\end{align*}
This formula will be used later. In addition, it can be easily shown that the Laplacian on $0$- and $1$-forms actually agrees with the usual componentwise Laplacian up to a sign.
\section{A Result in Integrability by Compensation}

Later, we shall make repeated use of the following compensation result:

\begin{Lm}
\label{prel:compensation}
	Let $da \in L^{m,\infty}(\R^{m})$, $db \in L^{p,r}(\R^{m})$ for $1 < p < + \infty$ and $1 \leq r \leq + \infty$. Then, we have $da \wedge db \in W^{-1, (p,r)}(\R^{m})$ together with the following estimate:
	\begin{equation}
	\label{compensationestimate}
		\| da \wedge db \|_{W^{-1, (p,r)}} \leq C \| da \|_{L^{m,\infty}} \| db \|_{L^{p,r}},
	\end{equation}
	for a constant $C > 0$.
\end{Lm}
\noindent
{\bf Proof of Lemma \ref{prel:compensation}.}

\noindent
	By density, we may assume $a,b \in \mathcal{S}(\R^{m})$, the general case follows by approximation. Let now $u$ be a solution of the following equation:
	\begin{equation}
		\Delta u = da \wedge db \text{ in } \mathcal{D}^\prime(\R^{m}).
	\end{equation}
	We will show that $\nabla u \in L^{p}$ as well as:
	\begin{equation}
		\| \nabla u \|_{L^{p}} \leq C \| da \|_{L^{m, \infty}} \| db \|_{L^{p}},
	\end{equation}
	the general case is a direct consequence of real interpolation (consider $da$ fixed to obtain the required linear operator in the interpolation argument). We distinguish two cases:
	\paragraph{Case 1:} If $p > \frac{m}{m-1}$, we know by the general H\"older inequality:
	\begin{equation}
		da \wedge db \in L^{q,r},
	\end{equation}
	where:
	$$\frac{1}{q} = \frac{1}{p} + \frac{1}{m}, \quad r = p.$$
	By elliptic regularity, we deduce that $u \in W^{2, (q,r)}$ and by Sobolev embeddings:
	$$\nabla u \in L^{p,p} = L^{p},$$
	together with the estimate:
	$$\| \nabla u \|_{L^{p}} \lesssim \| u \|_{W^{2, (q,r)}} \lesssim \| \Delta u \|_{L^{q,r}} \lesssim \| da \|_{L^{m, \infty}} \| db \|_{L^{p}}.$$
	
	\paragraph{Case 2:} If $p < \frac{m}{m-1}$, we take $\bar{b} \in \R$ such that $b - \bar{b} \in L^{p^\ast, p}$. Here, we denote by $p^\ast$ the parameter determined by:
	$$\frac{1}{p^\ast} = \frac{1}{p} - \frac{1}{m}.$$
	Observe that $da \wedge db = d\left( da \wedge (b - \bar{b}) \right)$. H\"older's inequality immediately shows:
	\begin{equation}
		da \wedge (b - \bar{b}) \in L^{q,p},
	\end{equation}
	where $1/q = 1/p^\ast + 1/m = 1/p$. Thus $q = p$. We therefore conclude:
	\begin{equation}
		\| da \wedge db \|_{W^{-1,p}} \lesssim \| da \wedge (b - \bar{b}) \|_{L^{p}} \lesssim \| da \|_{L^{m, \infty}} \| b - \bar{b} \|_{L^{p^\ast, p}} \lesssim \| da \|_{L^{m, \infty}} \| db \|_{L^{p}}.
	\end{equation}
	Thus, we may deduce:
	\begin{equation}
		\| \nabla u \|_{L^{p}} \lesssim \| da \|_{L^{m, \infty}} \| db \|_{L^{p}}.
	\end{equation}
	This finishes our proof. We emphasize that, in particular, the "critical" case $p = \frac{m}{m-1}$ is obtained by interpolation. \hfill $\Box$\\

 \end{document}